

\magnification =\magstep1
\baselineskip =14pt

\centerline {\bf ARRANGEMENTS OF HYPERPLANES AND VECTOR BUNDLES ON $P^n$.}
\vskip 1cm

\centerline {\bf I.Dolgachev and M.Kapranov\footnote{*}
{\rm  Research of both authors was partially supported
 by the National Science Foundation.}}
\vskip 1cm

\beginsection Introduction.

Let $X$ be a smooth algebraic variety and $D$ be a divisor with
normal crossing on $X$. The pair $(X,D)$ gives rise to a natural

sheaf $\Omega^1_X (\log D)$
 of differential 1-forms on $X$ with logarithmic poles on $D$.
For each point $x\in X$ the space of sections of this sheaf in a small
neighborhood of $x$ is generated over ${\cal O}_{X,x}$ by regular 1-forms
and by forms $d\log f_i$ where $f_i=0$ is a local equation of an irreducible
component of $D$ containing $X$. This sheaf (and its exterior powers)
 was originally introduced by Deligne [De]
to define a mixed Hodge structure on the open variety $X-D$.
An important feature of the sheaf $\Omega^1_X (\log D)$ is that it is locally
free i.e. can be regarded as a vector bundle on $X$.

\vskip .2cm

In this paper we concentrate on a very special case
when $X=P^n$ is a projective space and $D=H_1\cup ...\cup H_m$
 is a union of hyperplanes in general position.
It turns out that the corresponding vector bundles are quite interesting
from the geometric point of view. It was shown in an earlier paper [K]
of the second author that in this case
  $\Omega^1_{P^n}(\log D)$ defines an embedding of $P^n$ into the
Grassmann variety $G(n, m-1)$ whose image becomes, after the Pl\"ucker
embedding, a Veronese variety $V^{m-3}_n$ i.e. a variety projectively
isomorphic to the image of $P^n$ under the map given by the linear system
of all hypersurfaces of degree $m-3$. If the case when the hyperplanes
osculate
a rational normal curve in $P^n$ the bundle $\Omega^1_{P^n}(\log D)$
coincides with the secant bundle $E_n^m$ of Schwarzenberger [Schw1-2].
 The corresponding
Veronese variety consists in this case of chordal $n-1$ -{\rm dim}ensional
subspaces to a rational normal curve in $P^{m-2}$.

\vskip .2cm

The main result of this paper  (Theorem 7.2) asserts that in the case
$m\geq 2n+3$ the arrangement of $m$ hyperplanes ${\cal H} =
\{H_1,...,H_m\}$ can be uniquely reconstructed from the bundle
 $E({\cal H}) = \Omega^1_{P^n} (\log \bigcup H_i)$ unless all
its hyperplanes osculate the same rational normal curve of
 degree $n$.
To prove this we study the variety $C({\cal H})$ of jumping lines
 for $E({\cal H})$.
The consideration of this variety is traditional in the theory of
 vector bundles on $P^n$, see [Bar, Hu]. In our case this variety
is of some geometric interest.
For example, if $n=2$ i.e. we deal with $m$ lines in $P^2$ then
(in the case of odd $m$) $C({\cal H})$
is a curve in the dual $P^2$ containing the points corresponding to lines
from ${\cal H}$. The whole construction therefore gives a canonical way
to draw an algebraic curve through a collection of points in (the dual) $P^2$.

For 5 points $p_1,...,p_5$ in $P^2$ this construction gives the unique
conic through $p_i$.

 For 7 points $p_1,...,p_7$ the construction gives
 a plane sextic curve with for which $p_i$ are double points.
 A non-singular model of this curve is isomorphic to the plane
 quartic curve of genus 3 which is classically associated to 7
 points via Del Pezzo surfaces of degree 2, see [C2][DO].

For even number of lines in $P^2$ the set of jumping lines is
typically finite. In this case more interesting is the curve of
{\it jumping lines of second kind} introduced by Hulek [Hu].
 The study of this curve will be carried out in the
 subsequent paper [DK].

Let us only formulate the answer for 6 lines
(considered as points $p_1,...,p_6$ in the dual plane). In this case
Hulek's curve will be a sextic of genus 4 of which $p_i$ are nodes.
It is described as follows.
 Blow up the  points $p_i$.
 The result is isomorphic to a cubic surface  $S$ in $P^3$.
The inverse images if the points $p_i$  and
the strict preimages of quadrics through various 5-tuples of
$p_i$ form a Schl\"afli double sixer of lines on the cubic surface.
 To each such double sixer
there is classically associated a quadric $Q$ in $P^3$ called the
{\it Schur quadric} [Schur][R] (see also [B], p.162).
 It is uniquely characterized by the property
that the corresponding pairs of
lines of the double sixer are orthogonal with respect to $Q$. Our sextic curve
in $P^2$ lifted to the surface $S$ becomes the intersection $S\cap Q$.

In fact, much of Hulek's general theory of stable bundles on $P^2$ with
odd first Chern class can be neatly reformulated in terms of (suitably
generalized) Schur quadrics. This will be done in [DK].

\vskip .2cm

Thus our approach gives a unified treatment of many classical constructions
 associating a curve to a configuration of points in a projective space.
 It appears that a systematic study of logarithmic bundles in other situations
(like surfaces other than $P^2$)
will provide a rich supply of concrete examples and and give additional
insight into the geometry related to vector bundles.

We would like to thank L. Ein and H. Terao for useful
 discussions and correspondence related to this work.

\hfill\vfill\eject

\beginsection \S 1. Arrangements of hyperplanes.

\vskip 1cm

\noindent {\bf 1.1.} Let $V$ be a complex vector space of
 {\rm dim}ension $n+1$ and
$P^n = P(V)$ be the projective space of lines in $V$. Let ${\cal H} =
(H_1,...,H_m)$ be a set (arrangement) of hyperplanes in $P^n$.
Dually it defines a set (a configuration) of $m$ points in the
dual vector space $\check P^n = P(V^*)$. We say that ${\cal H}$
is in (linearly) general position if the intersection of any
$k\leq n+1$ hyperplanes from ${\cal H}$ is of co{\rm dim}ension
 exactly $k$. Throughout this paper we shall mostly deal with
arrangements in general position.

\vskip .3cm

\noindent {\bf 1.2.} We choose a linear equation $f_i \in V^*$
for each hyperplane $H_i$ of ${\cal H}$. This system of choices
defines a linear map
$$\alpha_{\cal H}: {\bf C}^m \rightarrow V^*;\quad
(\lambda_1,...,\lambda_m)\mapsto \sum \lambda_if_i.$$
The kernel of this map will be denoted by $I_{\cal H}$.
It consists of linear relations between the linear forms $f_i$.
By transposing, we obtain a linear map
$$\alpha_{\cal H}^*: ({\bf C}^m)^*
\longrightarrow I_{\cal H}^*.\eqno (1.1)$$
Assume that $m\geq n+2$ so that $I_{\cal H}\neq 0$.
 After making a natural
identification between the space ${\bf C}^m$ and its
 dual space $({\bf C}^m)^*$ (defined by the bilinear form
$\sum x_iy_i$) we obtain from the map (1.1) $m$
linear forms on the space $I_{\cal H}$ i.e. an arrangement
of hyperplanes in
$P({\cal H})$. We denote the arrangement thus obtained by
 ${\cal H}^{as}$ and refer to it as the {\it associated} arrangement.
It is clear that ${\cal H}$ is in general position if and only if the
restriction of the map $\alpha_{\cal H}$ to any coordinate subspace
 ${\bf C}^k$ in ${\bf C}^m$ with $k\leq n+1$ is of maximal rank.
 This implies that ${\cal H}^{as}$ is in general position if and only
if ${\cal H}$ is. In the latter case, the {\rm dim}ension of
$P(I_{\cal H})$ is equal to $m-n-2$. From now on whenever we
speak about the association we assume that the arrangements
are in general position. Note that to define the map
$\alpha_{\cal H}$ we need a choice of order on the set of
hyperplanes from ${\cal H}$. Making this choice we automatically
 make a choice on the set ${\cal H}^{as}$.

The notion of association was introduced by A.Coble [C1].
 For modern treatment see [DO]. This notion has been  rediscovered,
 under the names "duality" or
"orthogonality" several times later, notably in the context of
 combinatorial geometries (see [CR], \S 11) and hypergeometric
functions (see [GG]).

Obviously the association is a self-dual operation, so
$({\cal H}^{as})^{as} = {\cal H}$ where we make a canonical
identification between the spaces $V$ and $V^{**}$.

Let $V$ and $I$ be two vector spaces of {\rm dim}ensions $n+1$
and $m-n-1$ respectively and let ${\cal H}$ and ${\cal H}'$ be
two arrangements of $m$ hyperplanes
in $P(V)$ and $P(I)$ respectively. We shall say that ${\cal H}$
and ${\cal H}'$ are associated if there is a projective isomorphism
$P(I) \rightarrow P(I_{\cal H})$ taking ${\cal H}'$ to the associated
configuration ${\cal H}^{as}$ (matching the ordering, if it was made).
 In particular, when $m=2n+2$, we can speak about {\it self-associated}
 configurations.

By duality we can speak about associated configurations of points in
 projective spaces. The following proposition (equivalent to a result
 by A.Coble) gives a criterion of being associated (resp. self-associated)
in terms of the Segre (resp. Veronese) embedding. We state it in terms
 of configurations of points.

\proclaim 1.3.~Proposition. a) Let $V,I$ be vector spaces of
{\rm dim}ensions $n+1$ and $m-n-1$ respectively.
 Let $p_i\in P(V),\, q_i\in P(I), \,\, i=1,...,m$, be
 two configurations of $m$ points . Let $s(p_i,q_i)\in P(V\otimes I)$
 be the image
of  the pair $(p_i,q_i)$ with respect to the Segre embedding
 $s: P(V) \times P(I) \rightarrow P(V\otimes I)$.
The configurations of points $(p_1,...,p_m)$ and $(q_1,...,q_m)$
are associated to each other if and only if the points $s(p_i,q_i)$
 are projectively dependent but any proper subset of them is
 projectively independent. \hfill\break
b) Let $V$ be a vector space of {\rm dim}ension $n+1$ and
$p_i\in P(V),\, i=1,...,2n+2$ be a configuration of points.
Let $v(p_i) \in P(S^2V)$ be the image of $p_i$ under the
 Veronese embedding $v: P(V) \rightarrow P(S^2V)$.
 The configuration $(p_1,...,p_{2n+2})$ is self-associated
 if and only if the points $v(p_i)$ are projectively dependent
 but any proper subset of them is projectively independent.

\noindent {\sl Proof:} b) follows from a). For the proof of a),
see, e.g., [K].

\vskip .3cm

\noindent {\bf 1.4.} Let ${\cal H}$ be an arrangement of
 $m$ hyperplanes in $P(V)$ and ${\cal H}^{as}$ be the associated
arrangement in the space $P(I_{\cal H})$. Let
$$W = \{(\lambda_1,...,\lambda_m) \in {\bf C}^m:\,\,\sum\lambda_i =0\}.$$
Define a linear map
$$t_{\cal H}: I_{\cal H} \otimes V \rightarrow W$$
by the formula
$$t_{\cal H}( \, (a_1,...,a_m), v) = (a_1f_1(v), ... , a_mf_m(v)).$$
This map considered as an element of the tensor product
$I_{\cal H}^*\otimes V^*\otimes W$ will be of considerable
importance in the sequel. We shall refer to it as the
{\it fundamental tensor} of the configuration ${\cal H}$.

It is clear that the fundamental tensor
$t_{{\cal H}^{as}} \in (I_{{\cal H}^{as}})^*\otimes
 I_{\cal H}^*\otimes W = V^*\otimes I_{\cal H}^*\otimes W$
 of the associated configuration ${\cal H}^{as}$ is obtained
 from the fundamental tensor
$T_{\cal H}\in I_{\cal H}^*\otimes V^*\otimes W$ by interchanging
of factors in the tensor product.

In coordinates, fixing a basis $e_1,...,e_{n+1}$ in $V$ and its
 dual basis in $V^*$, let
 $A = \|a_{ij}\|_{1\leq i\leq n+1,\, 1\leq j\leq m}$
 be the matrix whose columns are the coordinates of
 the linear functions $f_i$ and
$B = \|b_{ij}\|_{1\leq i\leq m-n-1,\, 1\leq j\leq m}$
 be similar matrix for the associated arrangement.
 We can choose $B$ in such a way that $B\circ A = 0$.
 Then the coordinates of the tensor $t_{\cal H}$ are
 given by the formula
$$(t_{\cal H})_{ijk} = b_{ik}a_{kj}.$$

\proclaim 1.5.~Proposition. Suppose that ${\cal H}$ is in
general position. Then for any non-zero vector $v\in V$ the
linear operator $t_{\cal H}(v): I_{\cal H} \rightarrow W$
 defined by the fundamental tensor $t_{\cal H}$, is injective.

\noindent {\sl Proof:} If $(a_1,...,a_m)\in {\rm Ker} (t_{\cal H}(v))$
 then $a_if_i(v) = 0$ for all $i=1,...,m$. Let $J = \{i: f_i(v) = 0\}$.
Then for any $i\notin J$ we have $a_i=0$. Since ${\cal H}$ is
in general position, $|J|\leq n$. Hence
$\sum a_if_i=0$ is a non-trivial linear relation between $\leq n$
linear functions among $f_i$. This contradicts the assumption of
 general position for ${\cal H}$.

\vfill\eject

\beginsection \S 2. Logarithmic bundles.

\vskip 1cm

\noindent {\bf 2.1.}
Let ${\cal H}  = (H_1,...,H_m)$ be an arrangement of $m$
hyperplanes in $P^n = P(V)$ in general position.
We shall define the divisor $\bigcup H_i$ also by ${\cal H}$.
 This divisor has normal crossing. This means that for any
 point $x\in P^n$ its local equation can be given by $t_1...t_k=0$
where $t_1,...,t_k$ is a part of a system of local parameters at $x$.
 In this situation one can define the sheaf
$\Omega^1_{P^n} (\log {\cal H})$ of differential 1-forms with
logarithmic poles along ${\cal H}$, see [De]. It is a subsheaf
of the sheaf $j_*\Omega^1_U$ where $U = P^n -{\cal H}$
and $j:U\hookrightarrow P^n$ is the embedding.
If $x\in P^n$ and $t_1...t_k=0$ is a local equation of the
divisor ${\cal H}$ near $x$, as above, then the section of
$\Omega^1_{P^n} (\log {\cal H})$ near $x$ are meromorphic
differential forms which can be expressed as
$\omega + \sum u_i d\log t_i$ where $\omega$ is a 1-form
and $u_i$ are functions , all regular near $x$. It is not
 difficult to see that the sheaf $\Omega^1_{P^n} (\log {\cal H})$
 is locally free of rank $n$, see [De].

We shall denote the sheaf $\Omega^1_{P^n} (\log {\cal H})$ by
$E({\cal H})$ and call it the {\it logarithmic bundle} associated
 to ${\cal H}$. It will be the main object of study in this paper.
 We will not make a distinction between vector bundles and locally
free coherent sheaves of their sections.

\vskip .3cm

\noindent {\bf 2.2.} The sheaf $E({\cal H})^*$ dual to $E({\cal H})$
 has a nice interpretation in terms of vector fields. We say that a
 regular vector field $\partial$ defined in some open subset $U\i P^n$
 is tangent to ${\cal H}$ if for any $x\in U$ the vector $\partial (x)$
 lies in the intersection of the tangent hyperplanes at $x$ to all $H_i$
containing $x$ (in particular, $\partial (x) = 0$ if $x$ is a point of
 $n$ -tuple intersection). Such fields form a coherent subsheaf in the
 tangent sheaf $T_{P^n}$. It is easy to see by local calculations that
this sheaf is isomorphic to the dual sheaf $E({\cal H})^*$.

\proclaim 2.3.~Proposition.  Let $\epsilon_i :H_i \hookrightarrow P^n$
be the embedding map. We have the canonical exact sequence of sheaves
 on $P^n$
$$0\rightarrow \Omega ^1_{P^n} \rightarrow E({\cal H}) \buildrel
 res \over\rightarrow \bigoplus_{i=1}^m \epsilon_{i*} {\cal O}_{H_i}
 \rightarrow 0 \eqno (2.1)$$
where $res$ is the Poincar\'e residue morphism defined locally by the
 formula
$$a_1d\log t_1 + ... + a_k d\log t_k + b_{k+1}dt_{k+1} + ... + b_ndt_n
\longmapsto (a_1(x), ..., a_k(x), 0,...,0)$$
where $(t_1,...,t_n)$ is a system of local coordinates at $x$ such that
 $t_1...t_k=0$ is a local equation of the divisor ${\cal H}$ at $x$.

\noindent {\sl Proof:} See [De].

\vskip .2cm

The next two propositions follow simply from the above exact sequence.

\proclaim 2.4.~Proposition. The Chern polynomial $c(E({\cal H})) =
\sum c_i(E({\cal H}))t^i$ of the bundle $E({\cal H})$ is given by
$$c(E({\cal H})) = (1-ht)^{-m+n+1}$$
where $h$ is the class of a hyperplane in $P^n$.
  In particular,
the determinant $\bigwedge ^n E({\cal H})$ is isomorphic to the
line bundle ${\cal O}(m-n-1)$ on $P^n$.

\proclaim 2.5.~Proposition. a) The space $H^0(P^n, E({\cal H}))$
 has {\rm dim}ension $m-1$ and consists of forms
$$\sum_{i=1}^m \alpha_i d\log f_i = d\log ( \prod_{i=1}^m f_i^{\alpha_i}),
\quad \alpha_i\in {\bf C}, \sum\alpha_i = 0.$$
b) More generally, ${\rm dim} \, H^0 (E({\cal H})(k)) =
(n+1){k+n-1\choose n} - {k+n\choose n} + m{k+n-1\choose n-1}$.\hfill\break
c) $H^i (E({\cal H})(k)) = 0$ for $1\leq i\leq n-2$ and any $k\in {\bf Z}$.

Note that we can now identify the space $H^0(P^n, E({\cal H}))$
 with the space $W$ introduced in n. 1.4.

\vskip .3cm

\noindent {\bf 2.6.} The logarithmic bundles can be obtained from
the bundle
$\Omega^1_{P^n}$ by applying elementary transformations of vector
 bundles.
These transformations were introduced first in the case of vector
 bundles over curves by A.Tyurin [T] and their general definition
is due to Maruyama [M1-2]. Let us recall this concept.

Let $E$ be a rank $r$ vector bundle over a smooth algebraic variety
 $X$ and $Z\i X$ - a hypersurface. Denote by $i:Z\rightarrow X$ the
embedding. Suppose that we have chosen some quotient bundle $F$ of
 the restriction $i^*E$. Then we have a surjective map of sheaves
 $E\rightarrow i_*F$ on $X$. We define the coherent sheaf
${\rm Elm}^-_{Z,F}$ as the kernel of this surjection.
It is easy to see that it is locally free of rank $r$ i.e.,
 can also be regarded as a vector bundle.
This bundle is called the elementary transformation of $E$
 along $(Z,F)$.

Note than when $E$ is a line bundle and $F=i^*E$ then
${\rm Elm}^-_{Z,F}$ is just the twisted sheaf $E(-Z)$.

\vskip .3cm

\noindent {\bf 2.7.} The bundle $E$ can be reconstructed from
its elementary transformation by applying the "inverse"
elementary transformation ${\rm Elm}^+_{Z,F}$. In the
 situation of n.2.6, the definition of ${\rm Elm}^+$ is as follows.
 Let $E(Z)$ be the sheaf whose sections are sections of $E$ with
simple poles along $Z$.
Then ${\rm Elm}^+_{Z,F}(E)$ is a subsheaf of $E(Z)$ whose sections
after multiplying by the local equation of $Z$ belong to
 ${\rm Elm}^-_{Z,F}(E) = {\rm Ker}\{E\rightarrow i^*F\}$.
It is easy to see that
${\rm Elm}^+$ and ${\rm Elm}^-$ are mutually inverse operations.

\vskip .3cm

\noindent {\bf 2.8.} For any $1\leq i\leq m$ let ${\cal H}_{\leq i}$
 be the truncated arrangement $(H_1,...,H_i)$. By definition,
${\cal H}_{\leq 0} = \emptyset$ and $E({\cal H}_{\leq 0}) = \Omega^1_{P^n}$.
The residue exact sequence from n.2.4 induces the following exact sequence
$$0\rightarrow E({\cal H}_{\leq i-1})\rightarrow E({\cal H}_{\leq i})
 \rightarrow \epsilon_{i*}{\cal O}_{H_i} \rightarrow 0.$$
Passing to the dual exact sequence and using the adjunction formula
we find the following exact sequence
$$0\rightarrow E({\cal H}_{\leq i})^* \rightarrow
E({\cal H}_{\leq i-1})^* \rightarrow \epsilon_{i*}
 {\cal O}_{H_i}(1) \rightarrow 0.\eqno (2.2)$$
Thus, by definition, we obtain

\proclaim 2.9.~Proposition. For each $i\leq m$ we have isomorphisms
$$E({\cal H}_{\leq i-1}) \cong {\rm Elm}^-_{H_i,
\epsilon_{i*}{\cal O}_{H_i}} (E({\cal H}_{\leq i}) ),$$
$$E({\cal H}_{\leq i})^* \cong {\rm Elm}^-_{H_i,
\epsilon_{i*}{\cal O}_{H_i}(1)} (E({\cal H}_{\leq i-1})^* ).$$

It is the second isomorphism which will be useful for us in \S 5 later.

\proclaim 2.10.~Proposition. Assume $1\leq m\leq n+1$. Then
$$E({\cal H}) \cong ({\cal O}_{P^n})^{\oplus (m-1)} \,\,\oplus
\,\, {\cal O}_{P^n}(-1)^{n+1-m}.$$

\noindent {\sl Proof:} Since $m\leq n+1$, we can choose
homogeneous coordinates
  $x_1,...,x_{n+1}$   in
$P^n = P(V)$ such that the hyperplane $H_i,\, 1\leq i\leq m$,
is given by the equality
$x_i=0$.
By Serre's theorem [H] coherent sheaves on $P^n$ correspond to graded
${\bf C}[x_1,...,x_{n+1}]$ -modules (modulo finite-{\rm dim}ensional ones),
the correspondence being given by ${\cal F} \mapsto \bigoplus
H^0(P^n, {\cal F}(i))$.
We shall describe the module corresponding to $E({\cal H})$.
Denote this module
by $M'$. The ring ${\bf C}[x_1,...,x_{n+1}]$ will be denoted
shortly by $A$.

Denote by $\xi = \sum x_i \partial/\partial x_i$ the Euler
vector field on $V$.
By ${\rm Lie}_\xi$ and $i_\xi$ we shall denote the Lie
derivative along $\xi$ and the
contraction of 1-forms with $\xi$.

Let $\tilde{\cal H}\i V$ be the configuration of coordinate hyperplanes
$\{x_i=0\}, i=1,...,m$. Let $M$ be the space of all global sections of the
sheaf $\Omega^1_V(\log \tilde{\cal H})$ on $V$. It is a graded $A$ -
module; the graded component $M_r$ consists of forms $\omega$ such
that ${\rm Lie}_\xi \omega = r\cdot\omega$.

It is clear that the space
of sections $H^0(P^n, E({\cal H})(p))$ can be identified with the
subspace in $M_r$ consisting of forms $\omega$ such that $i_\xi\omega =0$.
Hence our module $M'$ corresponding to $E({\cal H})$ is the kernel of
 the homomorphism $M\rightarrow A$ given by $i_\xi$.

The graded $A$ -module $M$ is free: it is  isomorphic to
 $A^m \oplus A^{n-m}(-1)$ with the basis
$d\log x_1,...,d\log x_m, dx_{m+1}, ..., dx_{n+1}$, the first
$m$ elements being in degree 0, the remaining ones
 in degree 1. Since
$i_\xi( d\log x_i) = 1, \, i_\xi (dx_i) = x_i$, we find that $M'$
is the kernel
of the homomorphism
$$A^m \oplus A^{n-m}(-1) \rightarrow A,\quad (a_1,...,a_n)
 \mapsto \sum_{j=1}^m
a_j + \sum_{j=m+1}^{n+1} a_jx_j.\eqno (2.3)$$
However, an element $(a_1,...,a_n)$ from the kernel of (2.3)
 is uniquely determined by the components $(a_2,...,a_n)$ which
may be arbitrary: we just define $a_1$ to be equal
$-(\sum_{j=2}^m
a_j + \sum_{j=m+1}^{n+1} a_jx_j)$. This means that $M$ is
isomorphic to $A^{m-1} \oplus A(-1)^{n-m+1}$ as a graded $A$ -
module. Hence, by Serre's theorem,
$E({\cal H}) \cong ({\cal O}_{P^n})^{\oplus (m-1)} \,\,\oplus
 \,\, {\cal O}_{P^n}(-1)^{n+1-m}$.

\vskip .3cm

\noindent {\bf 2.11.}
The logarithmic bundles can be used to define a map from the projective
space to a Grassmannian with the image isomorphic to a Veronese variety.
Let us explain this in more detail.

By a Veronese variety we mean a subvariety in a projective space
$P^{{n+d\choose n}-1}$ which is projectively isomorphic to the image
 of the Veronese mapping
$$v_{n,d}: P^n = P(V) \rightarrow P^{{n+d\choose n}-1} =
 P(S^dV).\eqno (2.4)$$

Let $E$ be a vector space of {\rm dim}ension $n+d$ and
 $G(n,E)$ -- the Grassmannian of $n$ {\rm dim}ensional linear subspaces
in $E$.
We shall often identify it with the Grassmannian $D(d,E^*)$ of $d$ -
{\rm dim}ensional subspaces in the dual spaces $E^*$. Consider its
Pl\"ucker embedding
$$G(n,E) \hookrightarrow P(\bigwedge^n E) = P^{{n+d\choose n}-1}.\eqno
(2.5)$$
Note that the {\rm dim}ensions of the ambient spaces for the Pl\"ucker
 embedding
and the Veronese embedding coincide. Therefore it makes sense to
speak about $n$ -{\rm dim}ensional Veronese varieties in the
Grassmannian $G(n,E)$. The following result, proven in [K],
shows that the logarithmic bundle $E({\cal H})$ defines an
 embedding of $P^n$ into a Grassmannian whose image is a Veronese variety.

\proclaim 2.12.~Theorem.
Let ${\cal H}$ be an arrangement of $m\geq n+2$ hyperplanes
in $P^n$ in general position. Denote by
$W$ the space $H^0(P^n, E({\cal H})) \cong {\bf C}^{m-1}$.
For any point $x\in P^n$ consider the subspace of $W$ consisting
of all sections vanishing at $x$, and let  $\phi_{\cal H}(x)$ be
 the dual subspace of $W^*$. Then:\hfill\break
a) The {\rm dim}ension of $\phi_{\cal H}(x)$ equals $n$ for all
$x\in P^n$;\hfill\break
b) The correspondence $x\mapsto \phi_{\cal H}(x)$ is a regular
 embedding $\phi_{\cal H}: P^n \hookrightarrow G(n, W)^*$. \hfill\break
c) The image $\phi_{\cal H}(P^n)$ in $G(n, W^*)$ becomes, after
 the Pl\"ucker embedding $G(m-n-1,W) \i P(\bigwedge^{n} W^*)$,
a Veronese variety.

In particular, $E({\cal H})$ is the inverse image of the bundle
 ${\cal S}^*$ on $G(n,W^*)$ where ${\cal S}$ is the tautological
subbundle over $G(n,W^*)$.

\proclaim 2.13.~Corollary. Assume that $m=n+2$. Then $E({\cal H})
\cong T_{P^n}(-1)$ where $T_{P^n}$ is the tangent bundle of $P^n$.

\noindent {\sl Proof:} Since ${\rm dim} (W) = n+1$, the map
$\phi_{\cal H}$ defines an isomorphism $P^n = P(V)\rightarrow
 G(n,W^*) = G(1,W) = P(W)$. In this case the tautological subbundle
${\cal S}$ on $G(n,W^*)$ is isomorphic to $\Omega^1_{P(W)}(1)$.
 Hence $E({\cal H})$ is isomorphic to $T_{P^n}(-1)$.

\vskip .3cm

\noindent {\bf 2.14.} Let us call a rank $n$ vector bundle $E$ on
 $P^n$
{\it normalized} if $c_1(E)\in \{0,-1,...,-n+1\}$. If $E$ is
any rank $n$ vector bundle $E$ on $P^n$ and $c_1(E)= na+b$ where
$a\in {\bf Z}, b\in \{0,-1,...,-n+1\}$ then we denote by $E_{norm}$
the normalized bundle $E(-a)$.

In our case $c_1(E({\cal H})) = m-n-1$  so
the normalized bundle $E_{norm}({\cal H})$ has the form $E({\cal H})(-d+1)$,
where $m = 1 + nd + r, 0\leq r\leq n-1$. Its first Chern class equals $r-n$.
The case when the first Chern class of the normalized bundle is zero
 i.e. when $m=nd+1$, will play a special role for us since many results
 below rely on  a good theory of jumping lines for bundles with
 $c_1 =0$, see [Bar].

\hfill\vfill\eject

\beginsection \S 3. Steiner bundles.

\vskip 1cm

Vector bundles of logarithmic forms turn out to belong
to a more general class
of bundles remarkable for the existence of a very simple resolution.

\proclaim 3.1.~Definition. A vector bundle $E$ on $P^n = P(V)$
is called a Steiner bundle if $E$ admits a resolution of the form
$$0\rightarrow I\otimes {\cal O}_{P^n}(-1)\buildrel \tau \over
\longrightarrow W\otimes {\cal O}_{P^n} \rightarrow E \rightarrow
0\eqno (3.1)$$
where $I$ and $W$ are vector spaces identified with the corresponding
trivial vector bundles.

The bundles of this type were considered earlier by several people,
see [E][BS]. The name "Steiner bundles" will be explained later in
this section.

Note that applying the exact sequence of cohomology, we immediately
obtain
$$W\cong H^0(P^n, E),\quad I\cong H^0(P^n, E\otimes\Omega^1_{P^n}(1)).
\eqno (3.2)$$

\vskip .2cm

\proclaim 3.2.~Proposition.   A vector bundle $E$ be
  Steiner bundle if and only if  the cohomology groups
$H^q(P^n, E\otimes \Omega^{p}(p))$ vanish for all  $q>0$ and also for
 $q=0, p>1$. The resolution (3.1) is defined functorially in $E$.
More precisely, the tensor $t$ is the only non-trivial differential
 $d_1^{-1,0}:E_1^{-1,0}\rightarrow E_1^{0,0}$ of the Beilinson spectral
 sequence with the first term
$$  E_1^{pq} = H^q(P^n, F\otimes \Omega^{-p}(-p)) \otimes {\cal O}(p)\
eqno (2.3)$$
converging to $E$ in degree $0$ and to $0$ in degrees $\neq 0$.

The proposition follows easily from  considering the Beilinson
spectral sequence, see [E], Proposition 2.2.

\proclaim 3.3.~Corollary. The property of being a Steiner bundle is an
open property.

\vskip .2cm

\noindent {\bf 3.4.} A map  $\tau$ between sheaves
 $I\otimes {\cal O}_{P^n}$ and $W\otimes {\cal O}_{P^n}$, an in (2.1),
is uniquely determined by a tensor
$$ t\in {\rm Hom}\, (V, {\rm Hom}\, (I,W)) = V^*\otimes I^* \otimes W.
\eqno (3.3)$$
This tensor should be such that the map $\tau$ is fiberwise injective.

Thus we see that the fundamental tensor $t_{\cal H}$ of an arrangement
of hyperplanes in $P^n$ (see n.1.4) allows one to define a coherent
sheaf as the cokernel of the map
 $$\tau_{\cal H}: I_{\cal H} \otimes {\cal O}_{P^n}(-1) \rightarrow
W\otimes {\cal O}_{P^n}.\eqno (3.4).$$
Here the spaces $I = I_{\cal H}$ and $W$ are defined in nn. 1.2 and
1.4 respectively. It turns out that this sheaf is isomorphic to our
logarithmic bundle $E({\cal H})$.

\proclaim 3.5.~Theorem. Let ${\cal H}$ be an arrangement of $m$
hyperplanes in general position in $P(V)$. Suppose that $m\geq n+2$.
Then the logarithmic bundle $E({\cal H})$ is a Steiner bundle.
The corresponding tensor is the fundamental tensor $t_{\cal H}$ of
the configuration ${\cal H}$.

\noindent {\sl Proof:} Let $v\in V$ be a non-zero vector. The fiber
of the map
(3.4) over the point ${\bf C}v \in P(V)$  has, in the notation of
\S 1 the form
$$t_{\cal H}(v): I_{\cal H} \rightarrow W,\quad (a_1,...,a_m)
\mapsto (a_1 f_1(v),...,a_mf_m(v)).\eqno (3.5)$$
To prove our theorem, we shall construct, for any $v$, an explicit
isomorphism between ${\rm Coker}\, t_{\cal H}(v)$ and the fiber at
 ${\bf C}v$ of the bundle $E({\cal H})$. Consider the map of vector spaces
$$\pi_v:W \rightarrow E({\cal H})_{{\bf C}v},\quad (a_1,...,a_m)
\mapsto \sum a_i(d\log f_i) |_{{\bf C}v},\eqno (3.6)$$
where $E({\cal H})_{{\bf C}v}$ is the fiber of $E({\cal H})$ at
 ${\bf C}v$.
It follows from Theorem 2.12 a) that $\pi_v$ is a surjection.
Thus our theorem is a consequence of the following lemma.

\proclaim 3.6.~Lemma. We have ${\rm Ker}\, \pi_v = {\rm Im}\,
t_{\cal H}(v)$. In other words,
a section $\sum \lambda_i d\log f_i,\, \sum\lambda_i=0$ of
the bundle $E({\cal H})$ vanishes at ${\bf C}v$ if and only
if $\lambda_i = a_i f_i(v)$ for some
$(a_1,...,a_m) \in I_{\cal H}$.

\noindent {\sl Proof:} It suffices to show that ${\rm Im}\,
 t_{\cal H}(v) \i {\rm Ker}\, \pi_v$ since the spaces in
question have the same {\rm dim}ension.

Let $J = \{i: {\bf C}v \in H_i\}$ and $H_J = \bigcap_{i\in J}H_i$.
 A section $\omega$ of $E({\cal H})$ vanishes at ${\bf C}v$ if and
 only if $\omega$ is regular near ${\bf C}v$ as a 1-form and,
 moreover, vanishes on the tangent subspace to $H_J$.

Suppose that $\lambda_i =a_i f_i(v)$ where $(a_1,...,a_m)\in
 I_{\cal H}$. Then for $i\in J$ we have $\lambda_i = 0$ since
 $f_i=0$ on $H_i$. Hence
the form $\omega = \sum \lambda_i d\log f_i$ is regular at
 ${\bf C}v$. Let $\xi\in V$ be such that $f_i(\xi) =0$ for
 $i\in J$ i.e. $\xi$ represents a vector tangent to $H_J$ at
 ${\bf C}v$. Then the value of $\omega$ on this tangent vector
equals
$$\sum_{i\notin J} \lambda_i {f_i(\xi)\over f_i(v)} =
\sum_{i\notin J} a_i {f_i(v) f_i(\xi)\over f_i(v)} =
\sum_{i\notin J} a_i f_i(\xi) = \sum_{i=1}^m a_i f_i(\xi) = 0.$$
This proves the lemma and hence Theorem 3.5.

Let us mention that  is possible to give an alternative proof
of Theorem 3.5 by using Proposition 3.2.

\vskip .3cm

\noindent {\bf 3.7.} Let $E$ be a rank $r$ Steiner bundle on $P^n$
 given by the resolution (3.1). We have noticed already that
 $W = H^0(P^n,E)$. It is obvious that $E$ is generated by its
 global sections. Hence we obtain a regular map
$\gamma: P^n \rightarrow G(r, W^*)$ that takes a point $x\in P^n$
 into the dual of the subspace of sections vanishing at $x$.
This map can be defined
"synthetically" by means of the following "Grassmannian Steiner
 construction" [K].

Let $m={\rm dim} (W)+1$ so that ${\rm dim} (I) = m-1-r$. Take
$m-r-1$ projective subspaces $L_1,...,L_{m-r-1}$ in the projective
 space $P(W^*)$, each of co{\rm dim}ension $n+1$. Denote by $]L_i[$
 the "star" of $L_i$ i.e. the projective space of {\rm dim}ension
 $n$ formed by hyperplanes in $P(W^*)$ containing $L_i$. Identify
all the stars $]L_i[$ with each other by choosing projective
isomorphisms $\phi_i:P^n\rightarrow ]L_i[$. Suppose that for
any $x\in P^n$ the corresponding hyperplanes $\phi_i(x)$ are
 independent.
 Consider  the
locus of  subspaces in $P(W^*)$ of co{\rm dim}ension $m-r-1$
(i.e. of {\rm dim}ension $r-1$) which are intersections of the
corresponding hyperplanes from stars $]L_i[$ i.e. the subspaces of the form
$$\phi_1(x)\cap ...\cap \phi_{m-r-1}(x), x\in P^n.\eqno (3.7)$$
 This locus lies in $G(n,W^*)$.  It is   called the Grassmannian
 Steiner construction. This is s straightforward generalization
of the classical Steiner construction of rational normal curves,
 see [GH], Ch.4, \S 3. The following proposition shows
that this construction is equivalent to that of Steiner bundle.
This explains the name.

\proclaim 3.8.~Proposition. Let $X$ be a projective space of
{\rm dim}ension $n$ embedded in some way into the Grassmannian
$G^n(W)$ of co{\rm dim}ension $n$ subspaces in $W, \, {\rm dim} (W) = m-1$.
Let $Q$ be the rank $n$ bundle on $G^n(W)$, whose fiber over a subspace
$L\i W$ is $W/L$. Let $E$ be the restriction of $Q$ to $X$.
 The possibility of representing $X$ by the Grassmannian
 Steiner construction is equivalent to the fact that $E$ is a Steiner bundle.

\noindent {\sl Proof:}
The  choice of $m-n-1$ parametrized star $(]L_i[, \phi_i:P^n =
 P(V) \rightarrow ]L_i[)$
is equivalent to the choice of $m-n-1$ surjective linear operators
$a_i: W^*\rightarrow V^*$. Namely, given such  $a_i$, we associate
to any hyperplane in $V^*$ i.e. to any point of $P(V)$ its inverse
image under $a_i$. Thus a point $\pi\in G^n(W)$ corresponding to
$x\in P(V)$ is $\bigcap {\rm Ker} (a_i(x))$. Define a linear map
$A: {\bf C}^{m-n-1} \rightarrow {\rm Hom}\,(W^*,V^*)$ by setting
$a_i = A(e_i)$ where $e_1,...,e_{m-n-1}$ is the standard basis of
${\bf C}^{m-n-1}$. Denote the space ${\bf C}^{m-n-1}$ by $I$.
This defines a tensor $t\in I^*\otimes W\otimes V^*$ which, in
its turn, defines a morphism of sheaves
$I\otimes {\cal O}_{P(V)}(-1) \rightarrow W\otimes {\cal O}_{P(V)}$.
 Our bundle
$E$ must be the cokernel of this morphism. Indeed, dualizing, we have
to show that $E^*$ is the kernel of the dual map
 $W^*\otimes {\cal O}_{P(V)} \rightarrow I^*\otimes {\cal O}_{P(V)}(1)$.
 This is defined by a linear map
 $A^\dagger: V\rightarrow {\rm Hom}\,(W^*,I^*)$
 associated to $t$. The fiber of this bundle over a point
${\bf C}v$ of $P(V)$ is equal to the kernel of the linear map
$A^\dagger (v): L^*\rightarrow I^*$. The latter is dual to
 the point of $X\i G^n(W)$ corresponding to $x$. This identifies
 the fibres. The converse reasoning is obvious.

\proclaim 3.9.~Proposition {\rm [E]}. The rank of a non-trivial
Steiner bundle on $P^n$ is greater or equal to $n$.

\noindent {\sl Proof:} Let $r$ be the rank. A Steiner bundle
 is given by a linear map $V\rightarrow {\rm Hom}\,(I,W)$ where
${\rm dim} (I) = {\rm dim}(W) -r$. Let $D$ be the subvariety of
 ${\rm Hom}\,(I,W)$ consisting of linear maps of not maximal rank.
It is well known that its co{\rm dim}ension equals $r+1$
 (see [ACGH], p.67). Therefore if $V$ is of {\rm dim}ension $> r+1$
 every linear map $t:V\rightarrow {\rm Hom}\, (I,W)$ will map some non
- zero vector $v\in V$ to a matrix of not maximal rank. This does
not occur for Steiner bundles.

Thus logarithmic bundles provide examples of Steiner bundles
 of maximal possible rank. In the rest of this section we shall
 consider only rank $n$ Steiner bundles on $P^n$.

\vskip .3cm

\noindent {\bf 3.10.} Recall that a vector bundle $E$ is called
 stable if for any torsion - free coherent subsheaf $F\i E$ we have
$${\rm deg} (F) /{\rm rk} (F) < {\rm deg} (E) /{\rm rk} (E).$$
It is well known that the property of stability is preserved under
tensoring with invertible sheaves.

The following fact is a particular case of results of Bohnhorst and
Spindler ([BH], Theorem 2.7).

\proclaim 3.11.~ Theorem. Any non-trivial Steiner bundle on $P^n$ is stable.

\proclaim 3.12.~Proposition. Let $E$ be a non-trivial rank 3
Steiner bundle on $P^3$ with $c_1(E) = 3k$ (i.e. ${\rm dim} (W) = 3k+3$).
 Then the normalized bundle $E_{norm} = E(-k)$ is an instanton bundle on
$P^3$ i.e. $c_1(E_{norm}) =0$ and $H^1 (P^3, E_{norm}(-2)) =0$.

\noindent {\sl Proof:} Follows at once from the resolution (3.1).

\proclaim 3.15.~Corollary. Let ${\cal H}$ be a configuration of $m$
hyperplanes in $P^n$ in general position with $m>n+2$. Then:\hfill\break
a) The logarithmic bundle $E({\cal H})$ is stable.\hfill\break
b) If $n=3$ and $m=3d+1$ then the normalized bundle $E_{norm}({\cal H})
= E({\cal H}) (-d+1)$ is an instanton bundle on $P^3$.

\vskip .2cm

Denote by $M_{P^2}(a,b)$ the moduli space of stable rank 2 vector
bundles on $P^2$ with $c_1 =a, c_2 = b$. It is known to be an
irreducible algebraic variety
of {\rm dim}ension $4b-a^2-3$, see [OSS], Ch.2, \S 4. Note that
 $M_{P^2}(a,b)$ is isomorphic to the moduli space of normalized
 bundles namely to $M_{P^2}(0, (4b-a^2)/4)$ for $a$ even and to
 $M_{P^2}(-1, (4b-a^2+1)/4)$ if $a$ is odd.

\proclaim 3.16.~Corollary. If $a= m-3, \, b= {m-2\choose 2}$
for some $m$ then the moduli space $M_{P^2}(a,b)$ contains a
dense Zariski open subset consisting of Steiner bundles.

In other words, a generic stable bundle with these Chern classes
is a Steiner bundle.

\noindent {\sl Proof:} The property of being a Steiner bundle is
open (Corollary 3.3). The said moduli space indeed contains Steiner
bundles - the logarithmic bundles corresponding to configurations of
$m$ lines: they are stable by Theorem 3.11 and have the required Chern
classes by Proposition 2.2. Since the moduli space is irreducible, we are done.

Let us reformulate the above corollary in terms of normalized bundles.

\proclaim 3.17.~Corollary. For any $d>0$ each of the moduli space
 $M_{P^2}(0, d(d-1))$, $M_{P^2} (-1, (d-1)^2)$ has an open dense
subset consisting of twisted Steiner bundles.

\vskip .2cm

\noindent {\bf 3.18.} Notice that ${\rm dim} M_{P^2}(3,6) =
 {\rm dim} M_{P^2}(-1,4) = 12$. On the other hand, arrangements
of 6 lines in $P^2$ also depend on 12 parameters. We will show
 later that the map ${\cal H}\mapsto E({\cal H})$ from the space
 of arrangements of 6 lines to the moduli space ${\rm dim} M_{P^2}(3,6)$
 is generically injective. This will show that a generic bundle from
$ M_{P^2}(3,6)$ is a logarithmic bundle associated to an arrangement
of 6 lines in $P^2$.

\vskip .3cm

\noindent {\bf 3.19.} The operation of association discussed
in section 1 can be extended to Steiner bundles. Namely, we
can view the defining tensor (3.3) as a linear map $
V^*\otimes I^*\rightarrow W$  and consider the corresponding map
$$\tau ': V\otimes {\cal O}_{P(I)}(-1) \rightarrow W\otimes
{\cal O}_{P(I)}$$
of vector bundles on the projective space $P(I)$.

\proclaim 3.20.~Proposition - Definition. The map $\tau'$ is
injective on all the fibers if and only if $\tau$ is.
 In this case the Steiner bundle $ {\rm Coker} (\tau')$
is said to be associated to the Steiner bundle
$E = {\rm Coker} (\tau)$ and denoted by $E^{as}$.

\noindent {\sl Proof:} The condition that $\tau$ is
 not fiberwise injective
means that there are non-zero $v\in V,\, i\in I$ such
that $t(v\otimes i) = 0$. The same condition is equivalent
 to the fact that $\tau'$ is not fiberwise injective.
This proves the "proposition" part.

\vskip .2cm
The next proposition follows immediately from definitions
of n.1.4 and its proof is left to the reader.

\proclaim 3.21.~Proposition. Let ${\cal H}$ be an arrangement
of hyperplanes in $P(V)$ in general position and ${\cal H}^{as}$
be its associated arrangement in $P(I)$. Then there is a natural
isomorphism of vector bundles
$$E({\cal H}^{as}) \cong (E({\cal H}))^{as}.$$

\vfill\eject

\beginsection \S 4. Monoids, codependence and monoidal complexes.

\vskip 1cm

In this section we describe some constructions of projective geometry which
will be used in the study of logarithmic bundles, more precisely,
 in the description of jumping lines for such bundles.

\vskip .3cm

\noindent {\bf 4.1.} We shall work in projective space $P^n$ with
homogeneous coordinates $x_0,...,x_n$. Projective subspaces in $P^n$
 will be shortly called {\it flats}.
For a subset $S\i P^n$ let $<S>$ denote the flat (projective subspace)
spanned by $S$. In particular, for two points $p\neq q\in P^n$
 the notation $<p,q>$ means the line through $p$ and $q$.

Let $X$ be a hypersurface in a smooth algebraic variety $Y$ and
$x\in X$ be a point. We say that $x$ is a $k$ -tuple point of $Y$ i
f the whole $(k-1)$ -st infinitesimal neighborhood $x^{(k-1)}\i Y$
is contained in $X$.

As usual, if ${\cal L}$ is a line bundle on a projective variety $X$,
we shall denote by $|{\cal L}|$ the complete linear system of divisors on
 $X$ formed by zero loci of sections of ${\cal L}$ i.e.
 $|{\cal L}| = P(H^0(X,{\cal L}))$.

\vskip .3cm

\noindent {\bf 4.2.} Let $Z\i P^n$ be an irreducible variety.
 A hypersurface $X\i P^n$ of degree $d$ is called a $Z$ -monoid
 if each point of $Z$ is a
$(d-1)$ -tuple point of $X$. For example, a $Z$ -monoid of degree 2 is
just a quadric containing $Z$.

 We denote by $M_d(Z)$ the projective subspace of $|{\cal O}_{P^n}(d)|$
formed by all $Z$ -monoids of degree $d$.

We shall be mostly interested in the case when $Z\i P^n$ is a flat.
In this case, denoting $c={\rm codim}\, Z$, we find by easy
 calculation, that
$${\rm dim}\, M_d(Z) = {c+d-2\choose d-1} (n-c-1) +
{c+d-1\choose d} -1. \eqno (4.1)$$
In particular, if ${\rm codim}\, Z =2$ then ${\rm dim}\, M_d(Z) = nd$.

\proclaim 4.3.~Proposition. Let $Z\i P^n$ be a flat of
 {\rm dim}ension $k\leq n-2$. Any $Z$ -monoid is a rational
 variety ruled in $P^k$'s.

\noindent {\sl Proof:} Projecting $X$ from the subspace $Z$ we
find a rational map to $P^{n-k-1}$ whose fibers are flats of
 {\rm dim}ension $k$. Indeed, take any $(k+1)$ -{\rm dim}ensional
flat $L$ containing $Z$. Then $Z$ is a hyperplane in $L$.
 The intersection $L\cap X$ is a hypersurface of degree $d$
 in $L$ containing $d-1$ times the hyperplane $Z$. This means
that $L\cap X = (d-1)Z + H(L)$ where $H(L)$
is some hyperplane in $L$. So $H(L)$ is the fiber of the said
rational map over $L$, as claimed.

\vskip .3cm

\noindent {\bf 4.4.}
For any flat $Z\i P^n$  we denote by $]Z[$ the star of $Z$ i.e.
the projective space of hyperplanes containing $Z$, cf. n.3.8.
 Obviously ${\rm dim} \,\,]Z[ = {\rm codim}\, Z - 1$.

Assume that ${\rm codim} Z = 2$. There is a simple way to
construct irreducible $Z$ -monoids of degree $d$ by means
of the classical Steiner construction.
Take any point $x\in P^n - Z$ and any regular map
$$\psi: \,\,]Z[ \,\,\cong P^1 \rightarrow \,\,]x[\,\, \cong P^{n-1}$$
of degree $d-1$ ,i.e. a map given by a linear subsystem of
 $|{\cal O}_{P^1} (d-1)|$. Denote by $H$ the unique hyperplane
 containing $Z$ and $x$ and assume that $\psi (H) \neq H$.
 Consider the variety $X(Z,x,\psi)$ which is the union of
 codimension 2 flats $L\cap \psi (L),\,\, L\in ]Z[$. We claim
 that this is a $Z$ -monoid of degree $d$ containing the point $x$.

In fact, take a line $l$ which has no common points with
$Z \cup \{x\}$. Then
$]Z[$ is identified with $l$ by the correspondence taking
 $H\in \,\,]Z[$ to the intersection point $x_H = H\cap l$.
The map $\psi$defines a degree $d-1$ map $f: l\rightarrow l$
 defined as follows: $x_H \mapsto \psi (H)\cap l$. The graph
of this map $\Gamma_f \i l\times l \cong P^1\times P^1$
intersects the diagonal in $d$ points. Therefore
 $l$ intersects $X(Z,x,\psi)$ at $d$ points so
${\rm deg} \,\,X(Z,x,\psi) = d$. To see that $X(Z,x,\psi)$
  is a $Z$ -monoid, we take any general hyperplane $P$ containing
$Z$. The intersection $X(Z,x,\psi)\cap P$ is a hypersurface $Y$
in $P$ of degree $d$ which set-theoretically is the union of $Z$
and another hyperplane in $P$,the latter entering with multiplicity 1.
 This  implies that $Z$ enters with multiplicity $(d-1)$ into $Y$ so
$Z$ consists of $(d-1)$ -tuple points of $X(Z,x,\psi)$.

Conversely, every $Z$ -monoid containing a point $x$ outside $Z$
is equal to $Z(Z,x,\psi)$ for some regular map $\psi$ of degree $d$
from $]Z[$ to $]x[$. This follows from the proof of Proposition 4.3.
 So we have proving the following fact.

\vskip .2cm

\proclaim 4.5.~Proposition. Let $Z\i P^n$ be a codimension 2 flat
and $x\in P^n-Z$. Then there is a bijection between the set of
irreducible $Z$ -monoids of degree $d$ containing $x$ and the set
of regular maps $\psi: \,\,]Z[\,\,\rightarrow\,\, ]x[$ of degree
$d-1$ with the property that $\psi (<Z,x>) \neq <Z,x>$.

\vskip .3cm

\noindent {\bf 4.6.} We are going to relate monoids to a property
of point sets in projective spaces which we call {\it  codependence}.

Let ${\bf p} = (p_1,...,p_r)$ be an ordered $r$ -tuple of points
in $P^{n-1}$
and ${\bf q} = (q_1,...,q_r)$ be an ordered $r$ -tuple of points i
n $P^1$. We say that {\bf p} and {\bf q} are $d$ -{\it codependent}
 if there is a hypersurface $Y\i P^{n-1}\times P^1$ of bi-degree
$(1,d)$ which contains the points $(x_i,y_i)$. We say that {\bf p}
and {\bf q} are {\it strongly $d$ -codependent} if there is an
irreducible such hypersurface.

\vskip .2cm

\proclaim 4.7.~Proposition. Let $(\Lambda_1,...,\Lambda_r)$ be
an ordered
$r$ - tuple of hyperplanes in $P^{n-1}$ and $(q_1,...,q_r)$ be
an ordered $r$ -tuple of points in $P^1$. Let us regard each
$\Lambda_i$ as a point $p_i$ in the dual projective space
$\check P^{n-1}$. Then ${\bf p} = (p_1,...,p_r)$ and
${\bf q} = (q_1,...,q_r)$ are $d$ -codependent (resp.
strongly $d$ -codependent) if and only if there is a regular
map $\psi: P^1 \rightarrow P^{n-1}$ of degree $\leq d$ (resp.
of degree exactly $d$) such that $\psi (q_i) \in \Lambda_i$.

\noindent {\bf Proof:} Suppose that {\bf p} and {\bf q} are s
trongly $d$ -codependent. Let $Y$ be an irreducible hypersurface
of bi-degree $(1,d)$ in $\check P^{n-1} \times P^1$ containing
all $(p_i,q_i)$. Denote by $F(u,v) = F(u_0,...,u_{n-1}; v_0,v_1)$
 the bi-homogeneous equation of $Y$. We obtain, for any
$v = (v_0,v_1)$ a linear form $u\mapsto F(u,v)$ on $P^{n-1}$.
Since $Y$ is irreducible, this form is non-zero and so its kernel
is a hyperplane, denoted $\psi (v)$, in $\check P^{n-1}$ i.e. a
point in the initial $P^{n-1}$. This gives the desired map map
$\psi$ from $P^1$ to $P^{n-1}$. The rest of the proof is obvious.

\vskip .2cm

\proclaim 4.8.~Corollary. Let $Z\i P^n$ be a codimension 2
flat and $q,p_1,...,p_r$ be points in $P^n-Z$. The following
conditions are equivalent:
\item{(i)} There exists a $Z$ -monoid of degree $d$ (resp. an
irreducible $Z$ -monoid of degree $d$) containing $q,p_1,...,p_r$.
\item{(ii)} There exists a regular map $\psi:\,\, ]Z[\,\,
\rightarrow \,\,]q[$ of degree $\leq d-1$ (resp. of degree
exactly $d-1$) such that $\psi (<Z,p_i>)$ contains the line $<q,p_i>$.
\item{(iii)} The collection of points $<q,p_i> \in \,\,\,]q[
\,\,\cong P^{n-1}$ and $<Z,p_i> \in \,\,]Z[ \,\,\cong P^1$
are $(d-1)$ -codependent (resp. strongly $(d-1)$ -codependent).

The proof is immediate.

\proclaim 4.9.~ Proposition. For $r = nd$ all $(d-1)$ -
codependent pairs of $nd$ - tuples \hfill\break
$(x_1,...,x_{nd}, y_1,...,y_{nd})$  form a hypersurface $\Xi$
 in $(P^{n-1})^{nd}\times (P^1)^{nd}$. Let $x_{i0},...,x_{i,n-1}$
 be the homegeneous coordinates of the point $x_i\in P^{n-1}$ and
let $y_{i0},y_{i1}$ be the homogeneous coordinates of the point
$y_i\in P^1$. The equation of $\Xi$ is the
 determinant of size
$nd \times nd$ whose $i$ -th row is the following vector of length $nd$:
$$(x_{i0}y_{i0}^{d-1}, x_{i0}y_{i0}^{d-2}y_{i1},...,x_{i0}y_{i1}^{d-2};
 x_{i1}y_{i0}^{d-1}, x_{i1}y_{i0}^{d-2}y_{i1},...,x_{i_1}y_{i1}^{d-2};...$$
$$ ...;x_{i,n-1}y_{i0}^{d-1}, x_{i,n-1}y_{i0}^{d-2}y_{i1},...,
x_{i,n-1}y_{i1}^{d-1}).$$

\noindent {\sl Proof:} Let $V$ be the space of
polynomials in $(x_0,...,x_{n-1},y_1,y_2)$ homogeneous of degree 1
 in $x_i$ and of degree $d-1$ in $y_j$. Then $dim (V) = nd$.
The entries of the determinant in question are the values of
 $nd$ monomials forming a basis of $V$, on our $nd$ points.
 So the vanishing of the determinant is equivalent to the
linear dependence of the vectors given by these values i.e.
 to the $(d-1)$ -codependence of $(x_i)$ and $(y_i)$.

\vskip .3cm

\noindent {\bf 4.10.} Let $p_1,...,p_{nd+1}$ be $nd+1$ points in
 $P^n$ in general position. The {\it monoidal complex}
 $C(p_1,...,p_{nd+1})$ is, by definition the locus of all
 the codimension 2 flats $Z\i P^n$ for which there exists a
$Z$ -monoid of degree $d$ containing $p_1,...,p_{nd+1}$.

According to classical terminology of Pl\"ucker,
by complexes one meant 3-parametric families of lines in $P^3$
i.e. a hypersurface in the Grassmannian $G(2,4)$. As we shall see,
our $C(p_1,...,p_{nd+1})$ is a hypersurface in the Grassmannian
 $G(n-1, n+1)$. This explains the word complex.

Monoidal complexes will be used in the next section to
describe jumping lines of logarithmic vector bundles.
\vskip .2cm

\proclaim 4.11.~Theorem. Let $G = G(n-1, n+1)$ be the Grassmannian
of codimension 2 flats in $P^n$. Then:
\item{a)} For any points  $p_1,...,p_{nd+1}\in P^n$ in general
 position the variety $C( p_1, ..., p_{nd+1})$ is either the whole $G$ or
 a hypersurface in $G$. In the latter case its degree (i.e.  the
degree of any equation in Pl\"ucker coordinates defining
 this hypersurface) equals
$nd(d-1)/2$.
\item{b)} Any codimension 2 flat $Z$ containing one of the points
$p_i$ belongs to $C(p_1,...,p_{nd+1})$ and, moreover, is a $(d-1)$ -
tuple point of this variety.

The proof of this theorem will be organized as follows. As a first
step we
shall analyze the equation of $C( p_1, ..., p_{nd+1})$ and find its
 degree. Second step will be to prove part b) of the theorem, again
by using the equation.
These three steps will be done in nn. 4.12 and 4.13 respectively.

\vskip .3cm

\noindent {\bf 4.12.} We shall define codimension 2 flats by pairs
of linear forms whose coefficients are put into rows of a by
$2\times (n+1)$ matrix
$$ A = \pmatrix { a_{10}&a_{11}&...&a_{1n}\cr
a_{20}&a_{21}&...&a_{2n}}.$$
The flat corresponding to $A$ will be denoted $Z(A)$.
Its Pl\"ucker coordinates are just 2 by 2 minors of $A$.
A representation of a flat $Z$ as $Z(A)$ gives a parametrization
of the pencil
$]Z[$ of hyperplanes through $Z$ i.e. an explicit identification
$]Z[ = P^1$. Explicitly, to a point $(t_1,t_2)\in P^1$ we associate
 the hyperplane given by the equation
$$\sum_{j=0}^n (t_1a_{1j}+t_2a_{2j})x_j=0.$$
Denote the last point $p_{nd+1}$ by $q$.
 We can choose a coordinate system in $P^n$ in such a way that $q$
has coordinates $(1:0:...:0)$. The projective space $]{p_{nd+1}}[$
of lines through $q$
is identified with $P^{n-1}$. Explicitly, if $p=(b_0:...:b_n)\in P^{n}$
is another point then the line $<q,p>$ has homogeneous coordinates
$(b_1,...,b_n)$.

Let $b_{ij}, j=0,...,n$, be the homogeneous coordinates of the point
 $p_i\in P^n, i=1,...,nd$. The hyperplane $<Z(A),p_i>\in \,\,]{Z(A)}[$
has, under the above identification $]{Z(A)}[ \,\,=P^1$, the homogeneous
coordinates
$(\sum_j a_{1j}b_{ij}, \sum_j a_{2j}b_{ij})$.

Applying Corollary 4.8, we find that a flat $Z(A)$ belongs to the
variety $C(p_1,...,p_{nd+1})$ if and only if the two $nd$ -tuples
 of points
$$( (b_{i1},...,b_{in})\in P^{n-1}, i=1,...,nd)\quad {\rm and}
\quad ((\sum_j a_{1j}b_{ij}, \sum_j a_{2j}b_{ij})\in P^1, i=1,...,nd)$$
are $(d-1)$ -codependent. Substituting them into the determinant
of Proposition 4.9, we find an equation on matrix elements $a_{ij}$
whose degree in these elements equals $nd(d-1)$ (since each  entry
of the determinant will have degree $d-1$ in $a_{ij}$). The Pl\"ucker
 coordinates, being 2 by 2 minors of $A$, have degree 2 in $a_{ij}$.
Hence the degree of equation in Pl\"ucker coordinates equals
 $nd(d-1)/2$ as required.

\vskip .3cm

\noindent {\bf 4.13.} Let us prove part b) of Theorem 4.11.
 By symmetry it suffices to prove that if $p_{nd+1}\in Z$ then $Z$
 is a $(d-1)$ -tuple point of $C(p_1,...,p_{nd+1})$. Let us keep the
 notations and conventions introduced in the   proof of part a),
in particular, assume that $p_{nd+1} = (1:0;...:0)$. A flat $Z(A)$
 contains $p_{nd+1}$ if and only if $a_{10} = a_{20} =0$.

First, let us prove that such $Z(A)$ lies in $C(p_1,...,p_{nd+1})$.
 By Corollary 4.8, this means that the collections of points
$$(b_{i1}:...:b_{in})\in P^{n-1} \quad {\rm and} \quad (t_i,s_i) =
(\sum_{j=1}^n a_{1j}b_{ij},\, \, \sum_{j=1}^n a_{2j}b_{ij}) \, \in
\, P^1, \,\,\, i=1,...,nd$$
are always $(d-1)$ -codependent. To show this, we construct a polynomial
$F(b_1,...,b_n,t,s)$ homogeneous of degree 1 in $b_1,...,b_n$ and of degree
$(d-1)$ in $t,s$ such that for all $i$ we have $F(b_{i1},...,
b_{in},t_i,s_i)=0$ .
In fact, we can construct at least $d-1$ linearly independent such polynomials,
namely
$$F_m (b_1,...,b_n,t,s) = (\sum_{j=1}^n a_{2j}b_j) t^m s^{d-1-m} \,\, - \,\,
(\sum_{j=1}^n a_{1j}b_j)t^{m-1}s^{d-m},\quad m=1,2,...,d-1.$$
The possibility of finding $d-1$ such polynomials means that the kernel of
$nd \times nd$ -matrix whose determinant, by Proposition 4.9,
 defines $C(p_1,...,p_{nd+1})$, has dimension $\geq d-1$.
Since matrices with such properties are $(d-1)$ -tuple points
of the variety of degenerate matrices, we find that $Z(A)$ is
a $(d-1)$ - tuple point of $C(p_1,...,p_{nd+1})$.

So we have proven all the assertion of Theorem 4.11 except the
fact that the equation of the monoidal complex cannot be identically
 zero. This will be proven in \S 5.

\vskip .3cm

 Unfortunately, we do not know whether the case $C(p_1,...,p_{nd+1}) = G$
 occurs. It does not occur if $n=2$  (see Corollary 5.4). Another
case when this does not occur is as follows.

\proclaim Proposition 4.14. If $d\leq 3$ and $n$ is arbitrary then
 for any $nd+1$ points $p_1,...,p_{nd+1} \in P^n$ in linearly general
 position the monoidal complex $C(p_1,...,p_{nd+1})$ does not coincide
 with the whole Grassmannian $G$.

\noindent {\sl Proof:} Suppose the contrary i.e., that for any
codimension 2 flat $Z$ there is a $Z$ -monoid $X$ of degree $d$
 through $p_1,...,p_{nd+1}$. If we take $Z$ to lie in the hyperplane
$ H = <p_1,...,p_n>$ then we find that any such monoid $X$ should be
the union of $H$ and some $Z$ -monoid of degree $d-1$ through
$p_{n+1},...,p_{nd+1}$. Let $H' = <p_{n+1}, ..., p_{2n}>$.
We take $Z = H\cap H'$. By the above the corresponding monoid
should be the union of $H$, $H'$ and a $Z$ -monoid of degree
 $d-2$ through the remaining points $p_{2n+1},...,p_{nd+1}$.
If $d=2$ this means that $2n+1$ generic points $p_1,...,p_{2n+1}$
lie onthe union of two hyperplanes $H\cup H'$, which is impossible.
If $d=3$, the above means that $n+1$ points $p_{2n+1},..., p_{3n+1}$
 lie on a monoid of degree 1 i.e. on a hyperplane, which is also impossible.

 \vskip .3cm

\noindent {\bf 4.15. Examples in $P^2$.} Consider first the case $n=2$.
Then each $2d+1$ points in $P^2$ in general position define the
monoidal complex $C(p_1,...,p_{2d+1})$. It is a curve of degree
$d(d-1)$ with $(d-1)$ -tuple point at each $p_i$. Let us consider
 some particular cases.

a) Let $d=2$.  Then the curve $C(p_1,...,p_5)$ is just the
unique conic through the points $p_i$.

b) Let $d=3$. Then   $C(p_1,...,p_7)$ is a curve of degree 6 and genus 3
with double points at $p_i$. By definition, it is the locus of
 all possible singular points of  cubics through $p_1,...,p_7$.
This curve has the following (classical, see [DO]) description.

 Let $S\buildrel \sigma \over\rightarrow P^2$ be the blow-up of
$P^2$ in $p_1,...,p_7$. This is a Del Pezzo surface. Its anticanonical
linear system has dimension 2 and defines a double cover $S\buildrel
\tau \over\rightarrow \tilde P^2$ (this $\tilde P^2$ is different
 from the first one) ramified along a plane quartic curve $C'\i \tilde P^2$.
 We claim that $C'$ is birationally isomorphic to $C(p_1,...,p_7)$.

Indeed, the anticanonical linear system of $P^2$ consists of  cubic curves.
 The curves of the anticanonical linear system of the blown-up surface
 $S$ can be viewed, after projection $S\buildrel \sigma
 \over\rightarrow P^2$, as plane cubics through $p_1,...,p_7$.
 Denote this linear system by ${\cal L}\cong P^2$. The second
 projective plane $\tilde P^2$ is the space of lines in ${\cal L}$.
The projection $\tau$ associates to a point $p\in S$ with projection
 $z=\sigma (p)\in P^2$ the set of all plane cubics through
$p_1,...,p_7$ which also meet $p$ (so this set is a line in
 ${\cal L}$ i.e. a pencil of cubics). All the cubics from this
 pencil also contain some ninth point $p'$ which is conjugate
to $p$ with respect to the double cover $\tau$. The map $\tau$
 ramifies at $p$ when $p'=p$ i.e. the cubics from ${\cal L}$ have
a node at $z=\tau(p)$.

\vskip .2cm

c) Let $d=4$. Then  $C(p_1,...,p_9)$ is a curve of degree 12 with
triple points at $p_1,...,p_9$. Its genus equals 28. We de not know
 any special geometric significance of this curve.
\vskip .3cm

\noindent {\bf 4.16. Examples in $P^3$.} Let us consider the case $n=3$.
 The monoidal complex $C(p_1,...,p_{3d+1})$ is a line complex of degree
 $3d(d-1)/2$.

a) For $d=2$ i.e. for 7 points in $P^2$ we get the so-called
{\it Montesano complex}. It consists of lines in $P^3$ which
 lie on a quadric passing through points $p_1,...,p_7$. It is
 not difficult to see that this complex (which is a threefold
in $G(2,4)$) is isomorphic to a $P^1$ -bundle over a Del Pezzo
 surface of degree 2. This latter surface is obtained by blowing
 up the seven points of $P^2$ corresponding to $p_1,...,p_7 \in P^3$
 by association (see n.1.2 above). We refer to [Mo] for more details
about the geometry of this complex.

b) Let $d=3$. The complex $C(p_1,...,p_9)$ is a complex of degree 9
consisting of lines which appear as double lines of cubic surfaces
through $p_1,...,p_9$.

\vskip .3cm

\noindent {\bf 4.17. Examples in $P^n$.} We consider only the case
$d=2$ i.e. of $2n+1$ points in $P^n$. This case  gives a complex of
 codimension 2 flats which can be called the {\it generalized
Montesano complex}. We consider all quadrics in $P^n$ through
$p_1,...,p_{2n+1}$ and pick those among them which contain a
 $P^{n-2}$ i.e. quadrics of rank $\leq 4$. The union of all
$(n-2)$ -flats on all the quadrics of rank $\leq 4$ through
 $p_i$ gives a hypersurface in the Grassmannian $G(n-1,n+1)$
 whose degree equals $n$. This is our complex.

Note the case when all $p_i$ lie on a rational normal curve
$C$ in $P^n$ of degree $n$. In this case the generalized Montesano
complex will be the locus  of all $(n-2)$ -flats intersecting the curve $C$.
Its equation will be the Chow form of $C$ i.e. the resultant of two
 indeterminate polynomials of degree $n$. This is a consequence of
the following easy fact.

\proclaim 4.18.~Lemma. Let $C\i P^n$ be a rational normal curve and
$Z$ a codimension 2 flat in $P^n$. Then the two conditions are equivalent:
\item{(i)} There exists a quadric through $C$ and $Z$.
\item{(ii)} $C\cap Z \neq \emptyset$.

\noindent {\sl Proof:}  (ii)$\Rightarrow$(i):  Let $x\in C\cap Z$.
 Let $|{\cal O}(2)|$ be the linear system of all quadrics in $P^n$.
 We consider three projective subspaces $L, M, N \i |{\cal O}(2)|$
 consisting respectively of quadrics containing $C$,
containing $Z$ and containing $x$. Then $L,M\i N$. The
 codimension of $L$ in the whole $|{\cal O}(2)|$ is $2n+1$
 and hence its codimension in $N$ is
$2n$. The space $M$ has dimension $2n$. Hence $L\cap M \neq
\emptyset$ so there exists a quadric with required properties.

(i)$\Rightarrow$(ii): Let $Q$ be a quadric containing $C\cup Z$.
 Then $C$ and $Z$ are subvarieties in $Q$ of complementary dimensions.
In the case $n > 3$ (as well as in the case  when $n=3$ and $Q$ is
singular) this alone implies that the intersection is non-empty. If
$n=3$ and $Q$ is smooth the non-emptiness follows from the fact that
 $C$ regarded as a curve on $Q = P^1\times P^1$ has bidegree
(1,2) or (2,1). The case $n=2$ is trivial.

\hfill\vfill\eject

\beginsection \S 5. Monoidal complexes and splitting of logarithmic
bundles.

\vskip 1cm

\noindent {\bf 5.1.} One of the main tools for the study of vector
 bundles on $P^n$ is the restriction of bundles to projective subspaces
to $P^n$ especially to lines. By Grothendieck's theorem any vector
 bundle on $P^1$ splits into a direct sum of line bundles $\bigoplus
 {\cal O}_{P^1}(a_i)$.

 In this section we use this approach for logarithmic bundles
$E({\cal H})$
where ${\cal H} = (H_1, ... ,\hfill\break H_m)$ is an arrangement
 of $m$ hyperplanes in $P^n$ in general position. Let us write
  the number $m$ in the form
$m = nd+1+r$ where $d,r$ are integers and $0\leq r\leq n-1$.
Call a line $l$ in $P^n$ a {\it jumping line} for $E({\cal H})$
 (or for ${\cal H}$, if no confusion arises), if the restriction
 $E({\cal H})|_l$ is not isomorphic to ${\cal O}_l(d)^{\oplus r}
\oplus {\cal O}_l(d-1)^{\oplus (n-r)}$.

\vskip .2cm

Of special interest for us will be the case $m = nd  + 1$. In this
 case the normalized bundle $E({\cal H})_{norm} = E({\cal H})(-d+1)$
 has first Chern class 0. A line $l$ will be in this case a jumping
line for ${\cal H}$ if the restriction $E({\cal H})_{norm}|_l$ is non
-
trivial i.e. not isomorphic to ${\cal O}_l^{\oplus n}$.

\vskip .2cm

The main result of this section is as follows.

\proclaim 5.2.~Theorem. Suppose $m=nd+1$.  Let $p_1,...,p_{nd+1}$
be points in $\check P^n$ corresponding to hyperplanes
$H_1,...,H_{nd+1}\in H$. For any line $l\i P^n$ let $]l[$
be the corresponding codimension 2 flat in $\check P^n$.
Then a line $l\i P^n$ is a jumping line for the bundle
$E({\cal H})$ if and only if the flat $]l[$ belongs to the
 monoidal complex $C(p_1,...,p_{nd+1})$. In particular, the
locus of jumping lines of $E({\cal H})$ is the support of a
divisor in the Grassmannian $G(2,n+1)$ of degree $nd(d-1)/2$.

\proclaim 5.3.~Corollary. Assume $n=2$. Then the monoidal complex
 $C(p_1,...,p_{2d+1})$ (which is in this case a subvariety in the
 dual plane $\check P^2$), does not coincide with the whole $\check P^2$.

\noindent {\sl Proof:} This is a consequence of Theorems 5.2, 3.11
 and of the Grauert - M\"ulich theorem [OSS] which implies that the
 locus of jumping lines oif a stable rank 2 bundle on $P^2$ is in
 fact a curve.

\proclaim 5.4.~Corollary. Assume that $d\leq 3$. Then for any
configuration ${\cal H}$ of $nd+1$ hyperplanes
$H_1,...,H_{nd+1} \i P^n$ in general position the locus of jumping
 lines of the bundle $E_{norm}({\cal H})$ does not coincide with
the whole Grassmannian.

\noindent {\sl Proof:} This follows from Proposition 4.14.

Applying Corollary 4.8, we  can give an equivalent, more
 geometric description of the property of a line to be jumping.

\proclaim 5.5.~Corollary. Suppose $m=nd+1$. Let  $l\i P^n$  be a
 line intersecting the $H_i$ in distinct points. Then $l$ is a
jumping line for ${\cal H}$ if and only if there is a regular map
$\psi:l\rightarrow H_{nd+1}$ of degree $\leq d-1$ such that
 $\psi(l\cap H_i) \in H_{nd+1}\cap H_i$ for $i=1,...,nd.$

This reformulation is asymmetric: one of the hyperplanes,
namely $H_{nd+1}$, acts as a ``screen". Of course, any other
 $H_i$ can be chosen for this role.

\vskip .3cm

\noindent {\bf 5.6.} Let $E$ be a vector bundle on $P^n$.
We say that $E$ is projectively trivial if
 $E\cong {\cal O}_{P^1}(a)^{\oplus b}$ for some
$a\in {\bf Z},\, b\in {\bf Z}_+$. In this case the
 projective bundle $P(E)$ is trivial and, moreover, canonically
trivialized. To get the trivialization we note that $P(E) = P(E(-a))$.
If $W$ is the space of sections of $E(-a)$ then $E(-a)$ is canonically
identified with $W\otimes {\cal O}_{P^n}$. Hence $P(E)$ is canonically
 identified with $P^n \times P(W)$. For any two points $x,x'\in P^n$
we get the identification of fibres
$$\Psi_{E,x,x'}: P(E_x) \rightarrow P(W) \rightarrow P(E_{x'}).$$
We shall call this system of identifications the canonical
{\it canonical projective connection} of the projectively trivial bundle $E$.

\vskip .2cm

In our situation of logarithmic bundles it follows that whenever
$m=nd+1$ and $l$ is not a jumping line for ${\cal H}$, we get a
 canonical projective connection in the restricted bundle
$E({\cal H})|_l$. We are going to describe this connection explicitly.

Note that the fiber of the bundle $E({\cal H})^* =
 T_{P^n}(\log {\cal H})$ at any point $x\in P^n$ not
 lying on any $H_i$ is identified with the tangent space
 $T_xP^n$. Therefore the fibre  $P(E({\cal H}^*)_x)$ is
canonically identified with the projective space $P^{n-1}_x$
 of all lines through $x$. This means that for any non-jumping
 line $l$ the projective connection gives us isomorphisms which we denote
$$\Phi_{{\cal H}, l, x, x'}: P^{n-1}_x \rightarrow P^{n-1}_{x'},
\quad x,x' \in l -{\cal H}.$$
Our next result describes this identification.

\proclaim 5.7.~Proposition. Let $m=nd+1$ and let $l$ be a
 non-jumping line for ${\cal H}$. Let $x\in l-{\cal H}$ be
any point and $\lambda\in P^{n-1}_x$ be a line in $P^n$ through $x$.
 Then there is a unique regular map $\psi_{x,\lambda}: l\rightarrow
 H_{nd+1}$ of degree $d$ such that $\psi (l\cap H_i) \i H_{nd+1}
\cap H_i$ for each $i=1,...,nd$ and $\psi(x) = \lambda\cap H_{nd+1}$.
 For any other point $x\in l-{\cal H}$ the value at $\lambda$ of the
projective connection map $\Phi_{{\cal H}, l, x, x'}: P^{n-1}_x
\rightarrow P^{n-1}_{x'}$ equals the line
$<x',\psi_{x,\lambda}(x')> \in P^{n-1}_{x'}$.

\vskip .3cm

Now we start to prove our results. We shall begin with Theorem 5.2.
We need two lemmas.

\vskip .3cm

\proclaim 5.8.~Lemma.  A vector bundle $E^*$ on $P^1$ of rank $n$ and
first Chern class $(-n(d-1))$ does not have  the form
${\cal O}(-d+1)^{\oplus n}$ if and only if $H^0(P^1, E^*(d-2)) \neq 0$.

\noindent {\sl Proof:} As any bundle on $P^1$, our $E^*$
 has the form $\bigoplus _{i=0}^n {\cal O}_l(a_i)$ where
$\sum a_i = -n(d-1)$. The  condition $(a_1,...,a_n) \neq (-(d-1),...,-(d-1))$
is equivalent, under the above constaint on the sum, to
 the condition
"$\exists i : a_i\geq -d+2$" which is tantamount to
$H^0(P^1, E^* (d-2)) \neq 0$.
\vskip .2cm
The next lemma concerns the case when ${\cal H}$ consists
of just one hyperplane ${\cal H}$. In this case, as we have
seen in Proposition 2.10, the logarithmic bundles $E({\cal H})$
 is itself projectively trivial. So the canonical projective
connection on $E({\cal H})$ gives identifications:
$$\Phi_{H,x,x'}: P^{n-1}_x \rightarrow P^{n-1}_{x'},\quad x,x'\in P^n-H.$$
\proclaim 5.9.~Lemma. The identification $\Phi_{H,x,x'}$ takes
 a line $\lambda$ through $x$ to the line $<\lambda\cap H, x'>$ through $x'$.

\vbox to 5cm{}

\noindent {\sl Proof:} Let ${\bf H} \i {\bf C}^{n+1}$ be the linear
 hyperplane corresponding to the projective hyperplane $H\i P^n$.
 By Proposition 2.10, we have an isomorphism $E({\cal H}) \cong
 {\cal O}_{P^n} (-1)^{\oplus n}$. We can make this statement more
 precise by showing the existence of a natural isomorphism
$$E({\cal H}) \cong {\bf H}^* \otimes {\cal O}_{P^n} (-1).$$
Denote the space ${\bf C}^{n+1}$ shortly by $V$. Let $x$ be any
point of $P^n = P(V)$ and let ${\bf x}$ be the 1-dimensional subspace
 in $V$ representing $x$. The tangent space $T_xP^n$ is canonically
 identified with ${\bf x} \otimes V/{\bf x}$. Denote by $U$ the open
 set $P^n-H$. If $x\in U$ then the map ${\bf H} \rightarrow V
 \rightarrow V/{\bf x}$ is an isomorphism so we get identification
$T_xP^n = {\bf x}^* \otimes {\bf H}$. Correspondingly, the fiber at
 $x$ of $\Omega^1_{P^n}$ becomes identified with ${\bf x} \otimes
{\bf H}^*$ i.e., with the fiber at $x$ of ${\bf H}^* \otimes
 {\cal O}_{P^n} (-1)$. We get an isomorphism of restricted bundles
 $\phi: E({\cal H})|_U \cong {\bf H}^* \otimes {\cal O}_{P^n}(-1)$.
Using the fact that $E({\cal H})$ is isomorphic to
${\cal O}_{P^n}(-1)^{\oplus n}$, we can extend the isomorphism
$\phi$ to the whole $P^n$. In this model for $E({\cal H})$ the
fiber $P^{n-1}_x$ of $E({\cal H})^*$ at $x$ is canonically
identified with $H$ by assigning to the line $\lambda$ through
 $x$ the point $\lambda\cap H$. Our lemma follows from this immediately.

\vskip .3cm

\noindent {\bf 5.10.} Now we are ready to prove Theorem 5.2.
Let us consider the bundle $T_{P^n} (\log {\cal H})$ as the result
 of sucessive elementary transformations starting with the bundle
 \hfill\break
$T_{P^n} (\log \, H_{nd+1})$, as in Proposition 2.9. The latter
bundle is projectively trivial. Consider a line $l\i P^n$. We can
assume that $l\cap H_i$ are distinct points of $l$. Then the
restriction to $l$ of the bundle $T_{P^n}(\log\, {\cal H})$ is
 the elementary transformation of $T_{P^n} (\log\, H_{nd+1})  = {\cal
O}_l(1)^{\oplus n}$ with respect to points $y_i = l\cap H_i$ and subspaces
 $T_{y_i}H_i \i T_{y_i}P^n$.
\vskip .2cm

\proclaim 5.11.~Lemma. Consider on the projective line $P^1$ the
vector bundle ${\cal O}_{P^1}^{\oplus n} = {\cal O}_{P^1} \otimes E$
where $E$ is an $n$ -dimensional vector space. Let $y_1,..., y_{nd}$
be distinct points of $P^1$ and $\Lambda_1,...,\Lambda_{nd}$ be
 hyperplanes  in $E$. We regard $\Lambda_i$ as a hyperplane in the
fiber of our bundle over $y_i$. Then the following conditions are equivalent:
\item{(i)} The elementary transformation
 ${\rm Elm} _{\{y_1,...,y_{nd}\},\, \{\Lambda_1,...,\Lambda_{nd}\}}
(E\otimes {\cal O}_{P^1}(1))$ is not isomorphic to
 ${\cal O}_{P^1}(-d+1)^{\oplus n}$;
\item{(ii)} The two $nd$ -tuples $(\Lambda_1,...,\Lambda_{nd})
\in P(E)^{nd}$ and $(y_1,...,y_{nd}) \in (P^1)^{nd}$ are $d$ - \hfill\break
codependent in the sense of n. 4.6.

\noindent {\sl Proof:} Denote the elementary transformation in
condition (i) by
Elm. Then,
 By Lemma 5.8,, condition (i) is equivalent to nonvanishing of
  $H^0({\rm Elm}(E\otimes {\cal O}(1))(d-2))$. Let $x_0,x_1$ be
 homogeneous coordinates in $P^1$
A section of ${\rm Elm} (E\otimes {\cal O}(-1))(d-2)$ is, by
 definition, a homogeneous polynomial $s(x) = (x_0,x_1)$ of
degree $d-1$ with values in $E$ such that $s(y_i) \i \Lambda_i$.
 This is exactly the characterization of $(d-1)$ -codependence
 given in Proposition 4.7.

\vskip .2cm

\proclaim 5.12.~Corollary.  Let ${\cal H} = (H_1,...,H_{nd+1})$
 be a configuration of hyperplanes in $P^n$ in general position.
 A line $l\i P^n$  not lying in any $H_i$, is a jumping line for
${\cal H}$ if and only if the $nd$ -tuples
$(H_1\cap H_{nd+1},...,H_{nd}\cap H_{nd+1}) \in \check H_{nd+1}$ and
$(l\cap H_1,...,l\cap H_{nd+1})\in l^{nd}$ are $(d-1)$ -codependent.

\noindent {\sl Proof:} Let $l$ be given and suppose that $l$ does
 not lie in $H_{nd+1}$. Let ${\bf H}_{nd+1}$ be the linear hyperplane
 in ${\bf C}^{n+1}$ corresponding to $H_{nd+1}$.

By Proposition 5.7  the restriction
of $T_{P^n}(\log H_{nd+1})$ to $l$ is isomorphic to ${\bf H}_{nd+1}
\otimes {\cal O}_l(1)$.
The restriction of the bundle $E^*({\cal H})$ to $l$ is the
 elementary transformation of this projectively trivial bundle
with respect to points
$y_i = l\cap H_i$ and hyperplanes $\Lambda_i = T_{y_i}H_i$.
So we can apply Lemma 5.11. An explicit projective trivialization
of the
bundle $T_{P^n}(\log H_{nd+1})$ given in Proposition 5.7, identifies
the projectivizations of fiber at every point $y\in l-H_{nd+1}$ with
$H_{nd+1}$. Under this identification our hyperplanes $\Lambda_i$
 correspond to hyperplanes $H_{nd+1}\cap H_i \i H_{nd+1}$. So the
 assertion follows from Lemma 5.8

\vskip .3cm

\noindent {\bf 5.13.} To finish the proof of Theorem 5.2, let us
reformulate Corollary 5.12 in terms of the dual space $\check P^n$.
 Hyperplanes $H_i$ correspond to points $p_i$ of $\check P^n$,
the line $l$ corresponds to  a flat $Z$ of codimension 2.
The projective space $\check H_{nd+1}$ of hyperplanes in
 $H_{nd+1}$  becomes  the space $P^{n-1}_{p_{nd+1}}$ of lines
through $p_{nd+1}$ and $l$ itself becomes identified with the
 pencil $]Z[$ of lines through $Z$. Under these identifications
 the hyperplane $H_{nd+1}\cap H_i$ in $H_{nd+1}$ corresponds to the line
$<p_{nd+1}, p_i>$. Now Theorem 5.2 follows from the definition
 of the monoidal complex and Corollary 4.8.

\vskip .3cm

\noindent {\bf 5.14.} Proposition 5.7 now becomes just a reformulation
 of the fact that the projective connection in the projectively trivial
 bundle $E= {\cal O}_{P^1}(a)^{\oplus b}$ is induced by global sections
 of $E(-a)$. So it is proven.

\vskip .3cm

\noindent {\bf 5.15.} For any stable rank $r$ bundle $E$ on $P^n$ there
exists a Zariski open set $U \i G(2, n+1)$ such that for $l\in U$ the
splitting type $(a_1\geq ...\geq a_r)$ of $E|_U$ is constant.
A splitting type $(a_1,...,a_r)$ is called {\it generic} (or rigid)
 if $a_1 - a_r \leq 1$. By Grauert - M\"ulich theorem it is always
 generic if $r=2$. We conjecture that the splitting type of $E({\cal H})$
 is always generic.

\hfill\vfill\eject

\beginsection \S 6. Schwarzenberger bundles.

 \vskip 1cm

In this section we shall show that our logarithmic bundles
generalize the construction of Schwarzenberger [Schw1-2] of
vector bunldes of rank $n$ on $P^n$.

\vskip .3cm

\noindent {\bf 6.1.} Note that the following choices are equivalent:
\item{(i)} An isomorphism $ P^n \cong |{\cal O}_{P^1}(n)|$ of $P^n$
 with the $n$ -fold symmetric product of $P^1$.
\item{(ii)} A dual isomorphism $\check P^n \cong |{\cal O}_{P^1}(n)|^*$.
\item{(iii)} A map $\nu: P^1 \rightarrow \check P^n = |{\cal O}_{P^1}(n)|^*$
 given by a complete linear system.
\item{(iv)} A rational normal curve of degree $n$ (Veronese curve,
for short) $\check C$  in $P^n$, the image of the map in (iii).
\item{(v)} A map $\check \nu: P^1 \rightarrow P^n \cong |{\cal O}_{P^1}(n)|$
 given by a complete linear system.
\item{(vi)} A Veronese curve $C$  in $P^n$, the image of this map.

\vskip .2cm

\noindent Fix any of them. Then every point $x\in P^n$ is identified
 with a positive divisor $D_x$ of degree $n$ on $P^1$. Choose
 some $m\geq n+2$.  Let $V(x)$ be the subspace of sections
$s\in H^0(P^1, {\cal O}(m-2)$ whose divisor of zeroes ${\rm div} (s)$
 satisfies the condition ${\rm div} (s)\geq D_x$. Denote by
$V(x)^\bot \i H^0(P^1, {\cal O}(m-2))^*$ the orthogonal subspace.
Its dimension is equal to $n$. In this way we obtain a map
$$P^n \rightarrow G(n, H^0(P^1, {\cal O}(m-2))^*) = G(n,m-1),
\quad x\mapsto V(x)^\bot.\eqno (6.1)$$
Let $S$ be the tautological bundle on $G(n,m-1)$ whose fiber over
a point represented by an $n$ -dimensional linear subspace is this subspace.
The pull-back, with respect to (6.1), of $S$ is a rank $n$ vector
 bundle on $P^n$. It is defined by any of the above six choices,
in particular, by a choice of the Veronese curve $C\i P^n$.
We denote the dual bundle by $E(C,m)$ and call it the
{\it Schwarzenberger bundle} of degree $m$ associated to $C$.
 Thus fibers of $E(C,m)$ have the form
$$E(C,m)_x = { H^0(P^1, {\cal O}(m-2)) \over \{s\in H^0(P^1,
{\cal O}(m-2)): {\rm div} (s)\geq D_x\}}.\eqno (6.2)$$

If we fix another isomorphism $P^{m-2} \cong |{\cal O}_{P^1}(m-2)|$,
 this times by means of another Veronese curve $R$ of degree $m-2$
in $\check P^{m-2}$, then we can view each point $x\in P^n$ as a
 positive divisor $D_x$ on $R$ and the space $P(V(x)^\bot)$ as
 the projective subspace $<D_x>$ spanned by $D_x$ i.e. as an $(n-1)$ -
secant flat of $R$. For this reason the projective bundle $P(E(C,m))$
 is called the $n$ -{\it secant bundle} of $R$, see [Schw2].

One can easily show  that $E(C,m)$ is generated by its space of
global sections which is canonically isomorphic to
$H^0(P^1, {\cal O}(m-2)) = {\bf C}^{m-1}$.

\vskip .3cm

\noindent {\bf 6.2.} The bundle $E(C,m)$ is in fact a
Steiner bundle. This fact is well known, see e.g., [BS], Example 2.2.
Let us give a precise statement.

Fix an isomorphism $P^n = P(V), \, {\rm dim } \, V = n+1$. Then
 the choice of a Veronese curve $C$ is given by an isomorphism
 $V\cong S^n A$, where $A$ is a 2-dimensional vector space and
the points of the curve $C$ are represented by the $n$ -th powers
 $l^n, l\in A$. Consider the multiplication map
$$ t: V\otimes S^{m-n-2}A = S^n A \otimes S^{m-n-2}A \rightarrow
S^{m-2}A.\eqno (6.3)$$

\proclaim 6.3.~Proposition. The Schwarzenberger bundle $E(C,m)$ is
 a Steiner bundle on $P^n = P(S^nA)$ defined by vector spaces
$I = S^{m-n-2}A ,\, W = S^{m-2}A$ and the tensor $t: V\otimes I
 \rightarrow W$ given by (6.3).

\noindent {\sl Proof:} This follows from formula (6.2) and the fact
 that
for two sections $f\in H^0(P^1, {\cal O} (a))$, $\,\,$  $g\in H^0(P^1,
 {\cal O}(b))$ we have ${\rm div} (f) \geq {\rm div} (g)$ if and only
 if $f$ is divisible by $g$.
\vskip .2cm

\proclaim 6.4.~Theorem. Let ${\cal H} = (H_1,...,H_m)$ be an arrangement
 of $m$ hyperplanes in $P^n$ in general position. Suppose that all
 $H_i$ considered as points of $\check P^n$, lies on a Veronese curve
 $\check C$ (equivalently, all $H_i$ osculate the dual Veronese curve
 $C\i P^n$). Then there is an isomorphism
$$E({\cal H}) \cong E(C,m).$$

\noindent {\sl Proof:} This is equivalent to Theorem 3.8.5 from
[K] which describes the Veronese variety in the Grassmannian
corresponding to $E({\cal H})$. We prefer to give a direct proof here.

Let $I_{\cal H}$ be the space defined in n.1.4 and $W\i {\bf C}^{m}$
be the space of vectors with sum of coordinates zero, see formula (1.2).
 We shall construct explicit isomorphisms
$$\alpha: S^{m-n-2}(A) \rightarrow I_{\cal H},\quad \beta: S^{m-2}A
\rightarrow W $$
which take the multiplication tensor (6.3) into the fundamental tensor
 $t_{\cal H}$ which defines $E({\cal H})$ as a Steiner bundle.

\vskip .2cm

Let $f_i=0$ be the equation of the hyperplane $H_i$ from $H$.
 The condition that $H_i$ osculates $C$ means that $f_i$ considered
 as an element of $V^* = S^nA^*$ can be written in the form $u_i^n$
 where $u_i\in A^*$. Let $p_i\in P^1= P(A)$ be the point corresponding
 to $u_i$. Let us identify the space $S^{m-2}A = H^0(P^1, {\cal O}(m-2))$
 with the space $H^0(P^1, \Omega^1_{P^1}(p_1+ ...+p_m)$ of forms with
 simple poles at $(p_1,...,p_m)$. After that the map $\beta$ is
defined by the formula
$$\beta (\omega) = ({\rm res} _{p_1} (\omega), ...,
 {\rm res}_{p_m}(\omega)),\quad \omega \in H^0(P^1,
\Omega^1_{P^1}(p_1+ ...+p_m)).\eqno (6.4)$$
By the residue theorem the sum of components of $\beta
(\omega)$ equals 0 i.e.
$\beta (\omega)\in W$. Now, if we fix a point $q$ different
 from the $p_i$'s then we can identify the space
 $H^0(P^1, {\cal O}(n)) = S^n A$ with the space of rational
 functions on $P^1$ with poles of order $\leq n$ at $q$. Let
us denote this latter space by $L(nq)$. Let us also regard the
 space $S^{m-n-2}A$ as the space
$H^0(P^1, \Omega(p_1+ ...+ p_m -nq))$ of forms with at most
simple poles at $p_i$ and with zero of order $\leq n$ at $q$.
 This is a subspace of $H^0(P^1, \Omega^1_{P^1}(p_1+ ...+p_m)$
 and we define $\alpha$ to be the restriction of $\beta$ to this subspace.

Let us see that this is correct i.e.,the image of $\alpha$ indeed
lies in the space $I_{\cal H}$. By definition (see n.1.4),
$I_{\cal H}\i {\bf C}^m$ consists of $(\lambda_1,...,\lambda_m)$
 such that $\sum \lambda_if_i(v) =0$ for every $v\in V$. In
our case $\lambda_i = {\rm res}_{p_i}(\omega)$, where $\omega\in
H^0(P^1, \Omega(p_1+ ...+ p_m -nq))$. Since we have identified
 $V=S^nA = L(nq)$, we have, for any $v\in V$:
$$\sum\lambda_i f_i(v) = \sum \lambda_i (u_i^n, v) = \sum
{\rm res}_{p_i} (\omega)\cdot (u_i^n, v) = \sum {\rm res}_{p_i}
 (\omega\cdot v) = 0$$
since the form $\omega\cdot v$ has  poles only at $p_i$.

This shows the correctness of the definition of maps $\alpha$ and
 $\beta$. After the identification given by these maps it is obvious
that the multiplication tensor becomes identified with the fundamental
tensor $t_{\cal H}$.

\vskip .2cm

\proclaim 6.5.~Corollary. Suppose that ${\cal H, H}'$ are two
arrangements of $m$ hyperplanes in general position in $P^n$
such that all the hyperplanes from ${\cal H}$ and ${\cal H}'$
 osculate the some fixed Veronese curve $C\i P^n$. Then the
logarithmic bundles $E({\cal H})$ and $E({\cal H}')$ are isomorphic.

\proclaim 6.6.~Proposition. Let $C,C'$ be two Veronese curves
 in $P^n$ such that for some $m\geq n+2$ the Schwarzenberger
 bundles $E(C,m)$ and $E(C',m)$ are isomorphic. Then $C=C'$.

\noindent {\sl Proof:} We shall show that $C$ can be recovered
 intrinsically from $E(C,m)$.
Since $E(C,m)$ is a Steiner bundle, its defining tensor
 $t:V\otimes I \rightarrow W$ is determined by $E(C,m)$
 itself (see Proposition 3.2). We know that there are
isomorphisms $V\cong S^nA,\, I\cong S^{m-n-2}A$, $W \cong S^{m-2}A$,
 with ${\rm dim} A = 2$ which take the tensor $t$ into the multiplication
tensor (6.3). We shall see that as soon as such isomorphisms exist, the
 Veronese curves in $P(V), P(I), P(W)$ consisting of perfect powers
of elements of $A$, are defined by $t$ alone. Indeed, $t$ gives a
 morphism $T: P(V)\times P(I) \rightarrow P(W)$. The Veronese curve
 $R$ in $P(W)$ is recovered as the locus of $w\in P(W)$ such that
 $T^{-1}(w)$ consists of just one point, say $(v,i)$. The loci of
$v$ (resp. $i$) corresponding to various $w\in R$ constitute the
Veronese curves of perfect powers in $P(V)$ and $P(I)$ respectively.
 Proposition 6.6 is proven.

\vskip .3cm

\noindent {\bf 6.7.} Consider as an example the case when $m=n+3$.
It is well known that any $n+3$ points in general position in $P^n$
 lie on a unique Veronese curve,see [GH], p.530. So Corollary 6.5 and
 Proposition 6.6 lead to the following conclusion.

\item{} For two arrangements ${\cal H}$ and ${\cal H}'$ of $n+3$
hyperplanes in general position in $P^n$ the logarithmic bundles
$E({\cal H})$ and $E({\cal H}')$ are isomorphic if and only if the
Veronese curves osculated by ${\cal H}$ and ${\cal H}'$ coincide.

\noindent Moreover, any deformation of a bundle $E({\cal H})$ is
again of this type, as the following proposition shows.

\vskip .2cm

\proclaim 6.8.~Proposition. Any Steiner bundle $E$ on $P^n$ of
rank $n$ with ${\rm dim}\,\, I = 2,\, {\rm dim}\,\, W = n+2$ is
a Schwarzenberger bundle $E(C,n+3)$ for some Veronese curve $C\i P^n$.

\noindent {\sl Proof:} Let $E^{as}$ be the associated Steiner
bundle on $P^1$ (see n. 3.20). By Proposition 3.21, it suffices
to show that $E^{as} = E({\cal H})$ for some arrangement of points
on $P^1$. But this is obvious since $E^{as}$ is of rank 1 and hence
 is determined by its first Chern class, which is equal to $n+1$ in
 our case. Thus taking any $n+3$ points on $P^1$, we realize $E^{as}$
 as a logarithmic bundle and hence realize $E$ as a Schwarzenberger
 bundle.

\hfill\vfill\eject.

\beginsection \S 7. A Torelli theorem for logarithmic bundles.

\vskip 1cm

\noindent {\bf 7.1.} Let ${\cal A}_{gen}(m,n)$ be the variety of
all arrangements of $m$  \underbar {unordered} hyperplanes in $P^n$
in general position. So ${\cal A}_{gen}(m,n)$ is an open subset in
the symmetric product ${\rm Sym}^m(\check P^n)$ (Note that we do not
 factorize modulo projective transformations). The correspondence
${\cal H}\mapsto E({\cal H})$ defines a map
$${\cal A}_{gen}(m,n) \longmapsto M(n, (1-t)^{n+1-m}),\eqno (7.1)$$
where $M(n, (1-t)^{n+1-m})$ is the moduli space of stable rank $n$
 bundles on $P^n$ with Chern polynomial $(1-t)^{n+1-m}$. We are
interested in the question whether this map is an embedding. The
 statements that some moduli space is embedded into another are
 traditionally called "Torelli theorems" after the classical
Torelli theorem about the embeddine of the moduli space of curves
into the moduli space of Abelian varieties. The following theorem
 which is the main result of this section shows that $\Psi$ is very
 close to an embedding, at least for large $m$.

\vskip .2cm

\proclaim 7.2.~Theorem. Let $m\geq 2n+3$ and let ${\cal H,H}'$ be
two arrangements of $m$ hyperplanes in $P^n$ in general position.
Suppose that the corresponding logarithmic bundles $E({\cal H})$
and $E({\cal H}')$ are isomorphic. Then one of the two possibilities holds:
\item{1)} ${\cal H} = {\cal H}'$ (possibly after reordering the
 hyperplanes).
\item{2)} There exists a Veronese curve $C\i P^n$ such that all
hyperplanes from ${\cal H}$ and ${\cal H}'$ osculate this curve.
In this case $E({\cal H})$ and $E({\cal H}')$ are isomorphic to the
Schwarzenberger bundle $E(C,m)$.

\vskip .2cm

\noindent {\bf 7.3.} To prove Theorem 7.2 we have to recover
 (as far as possible) the configuration ${\cal H}$ from the bundle
$E({\cal H})$. The key idea is that the lines lying in each ${\cal H}_i$
are special jumping lines.

More precisely, we shall call a line $l\i P^n$ a {\it super-jumping line}
 for ${\cal H}$ if the restriction $E({\cal H})|_l$ contains as a
direct summand a sheaf ${\cal O}_l(a)$ with $a\leq 0$. As before,
let us write $m=nd+1+r$ with $0\leq r<n$. Clearly the line $l$ is
 super-jumping if and only if the restriction
$E_{norm}({\cal H})|_l$ of the  normalized bundle contains
${\cal O}_l(b)$ with $b\leq -d$.

\vskip .2cm

\proclaim 7.4.~Proposition. Any line $l$ lying in one of the hyperplanes
 $H_i$ of the arrangement ${\cal H}$ is a super-jumping line for ${\cal H}$.

\noindent {\sl Proof:} We can assume that $l\i H_1$. Let $F$ be a
 vector bundle on $P^1$. The property that $F$ contains
${\cal O}(a)$ with $a\leq 0$ as a direct summand, is equivalent
to the fact that $H^1(P^1, F(-2))\neq 0$.

We shall therefore prove that for $F=E({\cal H})|_l$ the above
cohomology does not vanish. Since the dimension of the cohomology
groups varies semi-continuously with $l$, it is enough for our
purpose to assume that $l\i H_1$ is not contained in any other
 $H_i, i\neq 1$. The residue exact sequence (2.1) gives a surjection
$$F\rightarrow {\cal O}_l \oplus \bigoplus_{i=2}^m {\bf C}_{l\cap H_i}
\rightarrow 0.$$
Hence $F(-2)$ maps surjectively onto ${\cal O}_l(-2)$. Since for
coherent sheaves on $P^1$ the functor $H^1$ is right exact, we get
a surjection $H^1 (F(-2)) \rightarrow H^1 ({\cal O}_l(-2)) = {\bf C}$.
 Proposition is proven.

\vskip .2cm

We want now to reformulate the condition of being a super-jumpimg
line in terms of the dual projective space $\check P^n$.

\vskip .3cm

\proclaim 7.5.~Proposition.  Let ${\cal H} = (H_1,...,H_{m})$
 be as before and let $l\i P^n$ be a line not lying in any $H_i$.
 Let $p_i$ be the point of the dual projective space  $\check P^n$
 corresponding to $H_i$. Let also $Z\i \check P^n$ be the codimension
2 flat corresponding to the line $l$. Then the following two conditions
 are equivalent:
\item{(i)} $l$ is a super-jumping line for ${\cal H}$;
\item{(ii)} There exists a quadric $Q\i P^n$ (of rank $\leq 4$)
containing all the points $p_1,..., p_{m}$ and the flat $Z$.

\noindent {\sl Proof:} Consider the dual bundle $E^*$ to
 $E=E({\cal H})$. In other words, $E^*$ is the bundle
 $T_{P^n}(\log {\cal H})$.  A line $l$ is super-jumping
for ${\cal H}$ if and only if
$H^0(l, E^*) \neq 0$. By reasoning analogous to those in
the proof of
Corollary 4.8 and Lemma 5.11, the existence of a section
of $E^*|_l$ is equivalent to the existence of a regular map
 $\psi: l\rightarrow H_m$ of degree 1 such that $\psi (l\cap H_i)
 \i H_m\cap H_i$ for $i=1,...,m-1$. A map $\psi$ is just an
identification of $l$ with some line $l'$ in $H_m$. Let $\Lambda,
 \Lambda'$ be codimension 2 flats in $\check P^n$ corresponding to
 $l,l'$. The map $\psi$ gives an identification $\Psi$ of the
projective lines (pencils) $]\Lambda[\,\, , \,\,\,]\Lambda'[\,$
 formed by hyperplanes through $\Lambda$ and $\Lambda'$ respectively.
 Such an identification defines, by Steiner's construction [GH], a
 quadric $Q$ of rank $\leq 4$. Explicitly, $Q$ is the union of
codimension 2 subspaces of the form $\Pi \cap \Psi (\Pi)$ where
 $\Pi\in\,\,\, ]\Lambda[$ is a hyperplane through $\Lambda$.
 This proves Proposition 7.5.

\vskip .3cm

\noindent {\bf 7.6.} We would like now to characterize the
hyperplanes $H_i$ of ${\cal H}$ as those of which every line
 is super-jumping for ${\cal H}$. To do this, it is again convenient
to use the dual projective space $\check P^n$ and the points
 $p_i\in \check P^n$ corresponding to $H_i$. Let us call a point
 $q\in \check P^n$ {\it adjoint} to $p_1,...,p_m$ if, $q$ does
 not coincide with any of $p_i$ and for any codimension 2 flat
$Z\i \check P^n$ containing $q$ there is a quadric containing $Z, p_1,...,p_m$.
For a fixed point $q\in \check P^n$ let $H\i P^n$ be the
 corresponding hyperplane. Proposition 7.5 shows that $q$ is
 adjoint to $p_1,...,p_m$ if and only if any line in $H$  is
 super-jumping for ${\cal H} = (H_1,...,H_m)$. Thus Theorem 7.2
is equivalent to the following fact (in which we write $P^n$
instead of $\check P^n$).

\vskip .2cm

\proclaim 7.7.~Theorem. Let $p_1,...,p_m$ be points in $P^n$
 in linearly general position and $m\geq 2n+3$. Then:
\item{a)} Unless all $p_i$ lie on one Veronese curve, there
are no points adjoint to $p_1,...,p_m$.
\item{b)} If all $p_i$ do lie on one Veronese curve $C$ then
points of $C$ and only they, are adjoint to $p_1,...,p_m$.

Part b) is a consequence of Lemma 4.21. We shall concentrate
on the proof of part b).
The proof will be based on the classical Castelnuovo lemma [GH].

\proclaim 7.8.~Lemma. Let $m\geq 2n+3$ and $p_1,...,p_{m}$ be
 points in $P^n$ in linearly general position which impose
 $\leq 2n+1$ conditions on quadrics. Then $p_i$ lie on a Veronese curve.

We shall prove the following fact which, together with Castelnuovo
 lemma will imply Theorem 7.7.

\proclaim 7.9.~Proposition. If $m\geq 2n$ and $q$ is adjoint to $p_1,...,p_m$
then $q,p_1,...,p_m$ impose exactly $2n+1$ conditions on quadrics.

\noindent {\sl Proof:} Let $L$ be the linear system of all quadrics through
first $2n$ points
$p_1,...,p_{2n}$. It is well known that $L$ has codimension $2n$ in
the linear system $|{\cal O}(2)|$ of all quadrics and quadrics from
$L$ cut out precisely $p_1,...,p_{2n}$. For any codimension 2 flat
$\i P^n$ let $M(Z)$ be the linear system of all quadrics through $Z$.
 It has dimension $2n$ (see n.4.2).
L
To establish Proposition 7.9 it suffices to prove the following fact.

\proclaim 7.10.~Proposition. Let $q$ be any point different from
 $p_1,...,p_{2n}$. Let $L_1\i L$ be the linear system of all quadrics
 through $p_1,...,p_{2n}, q$. Then:
\item{a)} For a generic codimension 2 flat $Z$ containing $q$ the
intersection $L\cap M(Z)$ consists of one point.
\item{b)} The points $L\cap M(Z)$ for generic $Z$ through $q$ as
above, span $L_1$ as a projective subspace.

\noindent {\sl Proof of 7.9 from 7.10:} Suppose we know Proposition 7.10.
Note that $L_1$ has codimension $2n+1$ in $|{\cal O}(2)|$. This follows
from the fact that quadrics from $L$ (in fact, just rank 2 quadrics from
 $L$) cut out
$p_1,...,p_{2n}$ and nothing else.

Let $L_2\i L_1$ be the linear system of quadrics through all
 $p_i,\,\,i=1,...,m$ and $q$. We shall show that $L_2 = L_1$.
 Indeed, suppose that $L_2$ is a proper
subspace in $L_1$. Since $q$ is adjoint to $p_1,...,p_m$, the
intersection $L_2\cap M(Z)$ is non-empty for any codimension 2
 flat $Z$ through $q$. But for generic such $Z$ the intersection
$L\cap Z = L_1\cap Z$ consists of just one point and these points
span $L_1$. Hence $L_2$ should miss some of these points, giving
a contradiction.

\vskip .3cm

\noindent {\sl Proof of Proposition 7.10 a):}  Suppose that the
statement is wrong. Then for every $Z$ through $q$ the linear
system $L\cap M(Z)$ contains a pencil. This means that for any
 additional point $r\in P^n$  and any $Z$ through $q$  there will
 be a quadric containing $p_1,...,p_{2n},r$ and $Z$. Now we shall
move $r$ and $Z$ in a special way to get a contradiction.

 Let us number the points $p_1,...,p_{2n}$ so
as to ensure that the hyperplane $<p_1,...,p_n>$ does not contain q.
Let $H = <q,p_1,...,p_{n-1}>$. Take  1-parameter families
 $Z=Z(t), r=
r(t)$, $t\in {\bf C}$ with the following properties:

\item{1)} $Z(t), r(t)$ lie in $H$ for any $t$.

\item{2)}  $Z(0) = <q, p_1,...,p_{n-2}>$,  the point $r(0)$ is a
generic point
inside $Z(0)$.

\item{3)} For $t\neq 0$ the flat $Z(t)$ intersects $Z(0)$ transversely
 inside $H$  and the curve $r(t)$ intersects Z(0) transversely at $t=0$.

Let $Q(t)$ be a quadric containing
$p_1,...,p_{2n}, r(t), Z(t)$ (which exists by our assumption). We can choose
$Q(t)$  for $t\neq 0$ to depend algebraically of $t$.
Let $Q = \lim_{t\rightarrow 0} Q(t)$. Then $Q$ contains the flat
 $<q, p_1,...,p_{n-2}>$ and, in addition,
the embedded tangent space to $Q$ at points $p_1,...,p_{n-2}$ and
$r(0)$
coincides the hyperplane  $H$. This means that the lines
$<p_i,p_{n-1}>$ and
$<r(0), p_{n-1}>$
will lie on $Q$. This implies that the whole hyperplane $H$ w
ill be part of $Q$.
Hence the remaining $n+1$ points $p_n,...,p_{2n}$ should lie on another
 hyperplane,
which is impossible. Part a) of Proposition 7.10 is proven.

\vskip .3cm

\noindent {\bf 7.11.} Now we shall concentrate on the proof of  part b)
of Proposition 7.10. Note that we can reformulate it as follows.

\proclaim 7.12.~Proposition. The linear system $L_1$ of quadrics through
$p_1,...,p_{2n}$ and $q$ is spanned by quadrics of rank $\leq 4$ contained
in this system.

Indeed, the union of $L_1\cap M(Z)$ for all codimension 2 flats $Z$ through
$q$
coincide with the part of $L_1$ consisting of quadrics of rank $\leq 4$.

So we shall prove Proposition 7.12. Note that is a consequence of the
 following fact.

\proclaim 7.13.~Lemma. The linear system $L$ of quadrics through
 $p_1,...,p_{2n}$ is spanned by $(1/2){2n\choose n}$ quadrics of
rank 2 contained in $L$.

Indeed, suppose we know Lemma 7.13. Let $Q$ be any quadric in $L_1$.
 By lemma,
$Q$ is a linear combination of rank 2 quadrics $Q_1,...,Q_N$ which
 lie in $L$ but not necessarily in $L_1$. Since $L_1$ is a hyperplane
in $L$, any pencil
$<Q_i,Q_j>$ intersects $L_1$ and $Q$ lies in the projective span
of their intersection points (for all $i,j$). But any quadric from
any pencil $<Q_i,Q_j>$
has rank $\leq 4$ since the $Q_i$ have rank 2. This reduces our
statement to Lemma 7.13.

\vskip .2cm

To prove Lemma 7.13,
let $D\i |{\cal O}(2)|$ be the locus of quadrics in $P^n$ of rank $\leq 4$.

\proclaim Lemma 7.14. The dimension of $D$ equals $2n$, the degree
of $D$ equals ${1\over 2}{2n\choose n}$ and $D$ spans the projective
 space $|{\cal O}(2)|$.

\noindent {\sl Proof of Lemma 7.14:} The space $D$ is the image of
 the map
$f: \check P^n \times \check P^n \rightarrow |{\cal O}(2)|$ which
takes $(H, H') \mapsto H\cup H'$. The map $f$ is generically two -
to - one. The inverse image under $f$ of
the standard bundle ${\cal O}(1)$ on the projective space
 $|{\cal O}(2)|$ is
${\cal O}(1,1)$ and its degree (self-intersection index) is
 ${2n\choose n}$. This
proves the statements about the dimension and the degree.
The fact that $D$ spans $|{\cal O}(2)|$ follows since $D$
 contains the Veronese variety of double planes which by
itself spans $|{\cal O}(2)|$.

\vskip .3cm

\noindent {\sl End of the proof of Lemma 7.13:}
 Note that $L$ intersects $D$ in finitely many points whose number
is equal to the degree of $D$. These points are, moreover, smooth
points of $D$ (being quadrics of rank exactly 2). Hence the i
ntersection is transversal at any of the points.

Suppose that the intersection points do not span $L$ and are
 contained in some hyperplane $M\i L$.

Take any quadric $X \in D$ which does not belong to $L$ and
consider the codimension $2n$ subspace $W$ in $L$ spanned by $M$ and $X$.
Then $W$ intersects $D$ in more points than the degree of $D$. This
means that the intersection will contain a curve (denote it $C(X)$) passing
through one of the $(1/2){2n \choose n}$
rank 2 quadrics which belong to $M$. Clearly there will be one of
 these quadrics, say, $Q$, which will be contained in any $C(X)$.

 The embedded tangent space $T_QC(X)$ is contained in the tangent
space to D at Q. Since $L$ intersects $D$ transversely, $T_QC(X)$ does
not intersect $M$ in any point other than $Q$. On the other hand, since
$D$ spans $|{\cal O}(2)|$, for a generic $X\in D$ the intersection of
the projective span $<M,X>$ with $T_QD$ will consist of $Q$ alone.
Since $T_QC(X) \i T_QD \cap \, <M,X>$, we get a contradiction.

This finishes the chain of reductions proving Proposition 7.10 b).
The proof of Theorem 7.7 is finished.

\hfill\vfill\eject

\beginsection References.

\vskip 1cm

\item{[ACGH]} E.Arbarello, M.Cornalba, P.Griffiths, J.Harris,
 Geometry of algebraic curves I (Grund. Math. Wiss, {\bf 267}),
 Springer-Verlag, 1985.

\item{[B]} H.F.Baker, Principles of Geometry, Vol. 3, Cambridge
Univ. Press, 1927.

\item{[Bar]} W.Barth, Moduli of vector bundles on the projective
plane, {\it Invent. Math.},
{\bf 42} (1977), 63-91.

\item{[BS]} G.Bohnhorst, H.Spindler, The stability of certain
vector bundles on $P^n$, in:
Lecture Notes in Math., {\bf 1507}, p.39-50, Springer-Verlag, 1992.

\item{[C1]} A.B.Coble, Associated sets of points, {\it Trans. AMS},
 {\bf 24} (1922), 1-20.

\item{[C2]} A.B.Coble,   Algebraic geometry and theta functions
 (AMS Colloquium Publ, Vol. 10), Providence RI 1929.

\item{[De]} P.Deligne, Th\'eorie de Hodge II, {\it Publ. Math. IHES},
{\bf 40} (1971), 5-58.

\item{[DK]} I.Dolgachev, M.Kapranov, Schur quadrics and their
generalizations in the theory of vector bundles on $P^2$, in preparation.

\item{[DO]} I.Dolgachev, D.Ortland, Point sets in projective
spaces and theta functions, {\it Asterisque} {\bf 165},
 Soc. Math. France, 1988.

\item{[E]} L.Ein, Normal sheaves of linear systems on curves,
in: {\it Contemporary Math.}, {\bf 116}, p.9-18, Amer. Math. Soc., 1991.

\item{[GG]} I.M.Gelfand, M.I.Graev, A duality theorem for
 general hypergeometric functions, {\it Sov. Math. Dokl.}
 {\bf 34} (1987), 9-13.

\item{[GH]} P.Griffiths, J.Harris, Principles of Algebraic geometry,
J.Wiley and Sons Publ., 1978.

\item{[H]} R.Hartshorne, Algebraic geometry, Springer-Verlag, 1977.

\item{[Hu]} K. Hulek, Stable rank 2 vector bundles on $P^2$ with
 $c_1$ odd, {\it Math. Ann.}, {\bf 242} (1979), 241-266.

\item {[K]} M.M.Kapranov, Chow quotients of Grassmannians I, to appear.

\item {[L]} J.Le Potier, Fibr\'es stables de rang 2 sur
 $P_2({\bf C}) $, {\it Math. Ann.} {\bf 241} (1979), 217-256.

\item {[M1]} M.Maruyama, Elementary transformations in the theory
of algebraic vector bundles, Lecture Notes in Math., {\bf 961},
 p.321-266, Springer-Verlag, 1983.

\item{[M2]} M.Maruyama, Singularities of the curve of jumping
lines of a vector bundle of rank 2 on $P^2$, Lecture Notes in Math.,
{\bf 1016}, p.370-411, Springer-Verlag, 1983.

\item{[Mo]} D.Montesano, Su di un complesso di rette di terzo grado,
{\it Mem. Accad. delle Scienze dell' Instituto di Bologna} (5),
 {\bf 3} (1893), 549 -577.

\item{[OSS]} C.Okonek, M.Schneider, H.Spindler, Vector bundles on
complex projective spaces, (Progress in Math., Vol.3), Birkh\"auser,
 Boston 1980.

\item{[R]} Th. Reye, Ueber die allgemeine Fl\"ache dritter Ordnung,
{\it Math. Ann.}, {\bf 55} (1901), 257-264.

\item{[Ro]} T.G. Room, The geometry of determinantal loci, Cambridge
Univ. Press, 1937.

\item {[Schur]} F. Schur, Ueber die durch collineare Grundgebilde
erzeugten Curven und Fl\"achen, {\it Math. Ann.}, {\bf 18} (1881), 1-32.

\item{[Schw1]} R.L.E.Schwarzenberger, Vector bundles on the projective
plane, {\it Proc. London Math. Soc.}, {\bf 11} (1961), 623-640.

\item{[Schw2]} R.L.E.Schwarzenberger, The secant bundle of a projective
variety, {\it Proc. London Math. Soc.}, {\bf 14} (1964), 369-384.

\item{[T]} A.N.Tyurin, On the classification of two - dimensional
 vector bundles over algebraic curves of arbitrary genus,
 {\it Izv. Akad. Nauk SSSR, ser. Math.}, {\bf 28} (1964), 21-52 (in Russian).

\vskip 2cm

Authors' addresses:

\noindent I.D.: Department of Mathematics, University of
Michigan, Ann Arbor MI 48109, email: IGOR.DOLGACHEV@um.cc.umich.edu

\vskip .3cm

\noindent M.K.: Department of Mathematics, Northwestern
 University, Evanston IL 60208, email: kapranov@chow.math.nwu.edu

\bye